\documentclass[12pt,leqno]{article}
\usepackage{lineno}
\usepackage{amsmath,amsthm}
\usepackage{txfonts}
\usepackage{bbm}
\usepackage{amssymb}
\usepackage{amsbsy}
\usepackage{mathrsfs}
\usepackage{graphicx}
 \usepackage{color,xcolor}

\allowdisplaybreaks
\def\ve{\varepsilon}

\theoremstyle{plain}
\theoremstyle{plain}
\theoremstyle{remark}  \newtheorem{remark}{\noindent\mbox{Remark}}
\theoremstyle{remark}  
\theoremstyle{plain}
\theoremstyle{plain}\newtheorem{lemma}{\noindent\mbox{Lemma}}
\theoremstyle{plain} \newtheorem{theorem}{\noindent\mbox{Theorem}}
\theoremstyle{plain}\newtheorem{proposition}{\noindent\mbox{Proposition}}
\theoremstyle{plain}\newtheorem{corollary}{\noindent\mbox{Corollary}}
\theoremstyle{definition} \newtheorem{definition}{\noindent\mbox{Definition}}
\theoremstyle{definition}
 \def\proof{\noindent{\it Proof.~~}}
 \def\qed{\hfill$\Box$\medskip}
 \def\rto{\rightarrow\infty}
 \def\z{\left}
 \def\y{\right}
 \def\no{\nonumber}

\begin{document}
 \title{{ \bf Beyond Poisson Approximation: Sums of Markovian Bernoulli Variables with Applications to Brownian Motions and Branching Processes}\footnote{Supported by Nature Science Foundation of Anhui Educational Committee, Anhui, China (Grant No.
2023AH040025).}}

\author{Hua-Ming \uppercase{Wang}$^{\dag,\ddag}$ and Shuxiong  \uppercase{Zhang}$^{\dag,\S}$  }
\date{}
\maketitle%
 \footnotetext[2]{School of Mathematics and Statistics, Anhui Normal University, Wuhu 241003, China}
\footnotetext[3]{Email: hmking@ahnu.edu.cn}
\footnotetext[4]{Email: shuxiong.zhang@mail.bnu.edu.cn}

\vspace{-.5cm}

\begin{center}
\begin{minipage}[c]{12cm}
\begin{center}\textbf{Abstract}\quad \end{center}

Let $\{\eta_i\}_{i\ge 1}$ be a sequence of dependent Bernoulli random variables. While the Poisson approximation for the distribution of $\sum_{i=1}^n\eta_i$
  has been extensively studied in the literature, this paper establishes new convergence regimes characterized by non-Poisson limits.  Specifically, under a Markovian dependence structure, we show that $\sum_{i=1}^n\eta_i,$ under suitable scaling, converges almost surely or in distribution as $n\to\infty$ to a geometric or Gamma random variable. These results provide a new tool for analyzing the limit distributions of sums of Markovian dependent Bernoulli random variables.
We demonstrate these results in several applications: determining the limiting distribution of the number of weak cutspheres for a $d(\ge3)$-dimensional standard Brownian motion; deriving the limit law for weak cutpoints of geometric Brownian motion; and analyzing how often the population size reaches a given threshold in certain branching processes, both with and without immigration.


\vspace{0.2cm}


\vspace{0.2cm}

\end{minipage}
\end{center}

 \section{Introduction}
 \subsection{Background, motivation and notation}
 The main purpose of this paper is to study
 the limit distribution of the number of successes in inhomogeneous Markovian Bernoulli trials. Let $\{\eta_i\}_{i\ge1}$ be a sequence of Bernoulli random variables. It is well-known that the Chen--Stein method, named after Chen \cite{c75} and Stein \cite{s72}, provides an efficient framework for bounding the approximation error between the distribution of
 $\sum_{i=1}^n\eta_i$ and a Poisson distribution. In this paper, from a different perspective, we aim to establish novel convergence regimes characterized by non-Poisson  limits. As applications, we study the number of weak cutspheres of a $d(\ge3)$-dimensional standard Brownian motion, the number of weak cutpoints of a geometric Brownian motion, and the number of times the population size of a Galton-Watson process or a branching process with immigration in varying environments reaches a given threshold.

  The approximation of $\sum_{i=1}^n\eta_i$ by a Poisson distribution has a long history  under various assumptions. For the case where $\{\eta_i\}_{i\ge1}$ are independent,  we refer the reader to \cite{hl60,k64,le60,pr53,v69} and references therein. Choi and Xia \cite{cc02} also showed that  the binomial approximation is generally better than the Poisson approximation for independent Bernoulli trials.  When $\{\eta_i\}_{i\ge1}$ are dependent, using an approach similar to that of Stein \cite{s72}, who studied normal approximation for dependent random variables, Chen \cite{c75} established a general method for bounding the error in the Poisson approximation of $\sum_{i=1}^n\eta_i.$  This methodology, now widely recognized as the Chen--Stein method, has proven to be remarkably effective for Poisson approximation problems. For further advancements in this area, see \cite{agg89,cr13,pp02,xia97} and references therein.

   It is important to note that the primary aim of these works is to bound the approximation error between the distribution of $\sum_{i=1}^n\eta_i$ (for fixed $n$) and a Poisson distribution.
For instance, in the case where $\{\eta_i\}_{i\ge1}$ are independent,  $P(\eta_i=1)=p_i$ and $\max_{1\le i\le n}p_i\le 1/4,$   Le Cam \cite{le60} proved that for any real-valued function $h$ defined on nonnegative integers, $
 \z|E\z(h\z(\sum_{i=1}^n\eta_i\y)\y)-\mathscr P_{\lambda}h\y|<16\lambda^{-1}\sum_{i=1}^np_i^2,$
where  $\mathscr P_\lambda h=\sum_{i=1}^\infty e^{-\lambda}\lambda^ih(i)/i!$ and $\lambda=\sum_{i=1}^np_i.$

In contrast, our work adopts a distinct perspective. When the sequence $\{\eta_i\}_{i\ge1}$ is governed by a Markovian dependence structure, we show that under suitable scaling conditions, the sum $\sum_{i=1}^n\eta_i$ converges almost surely or in distribution as $n\rto$ to a random variable with a geometric or gamma distribution.

Prior to presenting the main results, we introduce some notation and conventions.

\noindent{\bf Notation:} In this paper, we use $\xrightarrow{d}$ and $\xrightarrow{a.s.}$ to denote convergence in distribution and  almost sure convergence, respectively. We denote by $\mathbb N,$ $\mathbb R,$ and $\mathbb C$ the sets of natural, real and complex numbers, respectively. For $x\in \mathbb R,$ $\lfloor x\rfloor$ denotes the floor function, i.e., the greatest integer less than or equal to $x$. Furthermore, the asymptotic notation $f(n)\sim g(n)$ indicates that  $\lim_{n\rto}f(n)/g(n)=1.$ We write $f(n)=O(g(n))$ if there is a constant $c>0$ such that $|f(n)|\le cg(n)$ for all sufficiently large $n.$ For $\alpha>0,$ the Gamma function is defined as $\Gamma(\alpha):=\int_{0}^\infty x^{\alpha-1}e^{-x}dx.$  For $\beta,\ r>0,$  we write $\xi\sim\mathrm{Gamma}(r,\beta)$ if the random variable $\xi$ has the density function $\beta^{r}x^{{r-1}}e^{-\beta x}/\Gamma({r}),~x>0.$  Similarly, for $\lambda>0$ and $p\in (0,1),$ we write $\xi\sim\mathrm{Exp}(\lambda)$ if $\xi$ has the density function $\lambda{e}^{-\lambda x},\ x>0,$ and $\xi\sim\mathrm{Geo}(p)$ if $P(\xi=i)=(1-p)^{i}p$ for $i\geq 0$. Finally, for a set $A$, we denote by $\#A$ its cardinality.

\subsection{Main results}

To state the main results, we need to define the so-called regularly varying sequence.
\begin{definition} Let $\{S(n)\}_{n\ge 1}$ be a sequence of positive numbers. If for each $x>0,$
\begin{align*}
  \lim_{n\rto}\frac{S(\lfloor xn\rfloor)}{S(n)}=x^{\tau},
\end{align*}
we say the sequence $\{S(n)\}_{n\ge 1}$ is regularly varying with index $\tau.$

\end{definition}

Suppose that $\{\eta_n\}_{n\ge 1}$ is a sequence of Bernoulli variables such that  for $i\ge1,$ $k\ge 1,$ $1\le j_1\le\ldots\le j_k\le j_{k+1}<\infty,$ \begin{equation}\label{mc}\begin{split}
&P\z(\eta_i=1\y)=\frac{1}{\rho(0,i)},\ P\left(\eta_{j_{k+1}}=1\middle|\, \eta_{j_1}=1,\ldots,\eta_{j_k}=1\right)=P\left(\eta_{j_{k+1}}=1\middle|\,\eta_{j_k}=1\right)=\frac{1}{\rho(j_k,j_{k+1})}, \end{split}
\end{equation} where $\rho(j,j)=1,$ $\rho(i,j)> 1$ for $j>i\ge 0.$

We have the following theorem, which says that under certain conditions, when properly scaled,  $\sum_{j=1}^n\eta_j$ converges to a non-degenerate random variable as $n\rto.$

\begin{theorem}\label{thbb} Let $\{\eta_{n}\}_{\ n\ge 1}$ be a sequence of Bernoulli variables such that \eqref{mc} is satisfied. Assume that $\{D(n)\}_{n\ge0}$ is a sequence of numbers such that $D(0)=1$ and $D(n)>1$ for all $n\ge1.$
\begin{itemize}
  \item[{\rm(i)}]Assume that $\zeta(D):=\sum_{n=1}^\infty \frac{1}{D(n)}<\infty.$ If $\rho(i,j)=D(j-i)$ for all $j\ge i\ge 1,$  then
  $$\sum_{j=1}^n\eta_j\xrightarrow{a.s.}\xi \text{ as }n\rto,$$
  where $\xi\sim\mathrm{Geo}\Big(\frac{1}{\zeta(D)+1}\Big).$

  \item[{\rm(ii)}]
  Suppose that $\sum_{n=1}^\infty \frac{1}{D(n)}=\infty.$
       For $n\ge 1,$ let $S(n)=\sum_{i=1}^n\frac{1}{D(i)}.$ Assume $\{S(n)\}_{n\ge1}$ is regularly varying with index $\sigma\in [0,1].$ If for each $\ve>0,$ there exists $n_1\ge 1$ such that
   \begin{align}\label{sigma01} 1-\ve \le \frac{\lambda_\sigma D(i)}{\rho(0,i)}\le 1+\ve, \quad
     1-\ve\le \frac{D(j-i)}{\rho(i,j)}\le 1+\ve
  \end{align}
  for $j-i\ge n_1,\ i\ge n_1$ with $\lambda_0=1$ and $\lambda_\sigma=\frac{\Gamma(1+2\sigma)}{\sigma\Gamma(\sigma)\Gamma(1+\sigma)}$ for $\sigma\in (0,1],$  then
   $$\frac{\sum_{j=1}^n\eta_j}{S(n)}\xrightarrow{d}\xi \text{ as }n\rto,$$
    where $\xi\sim\mathrm{Exp}(\lambda_\sigma).$
    \end{itemize}
\end{theorem}
\begin{remark}
For $j\ge i\ge1,$  $\rho(i,j)$ characterizes the dependence between $\eta_i$ and $\eta_j.$ From \eqref{sigma01},  we note that to apply the theorem, it is not necessary to strictly require  $\rho(i,j)$ to be an explicit function of the distance $j-i.$ In practice, if $\rho(i,j)$ is asymptotically characterized by the distance between $i$ and $j,$ the theorem remains valid.
  \end{remark}
In certain special cases, even if condition \eqref{sigma01} fails, it is still possible to determine the limit distribution of $\sum_{i=1}^n \eta_i$.

\begin{theorem}\label{tha} Let $\{\eta_{n}\}_{\ n\ge 1}$ be a sequence of Bernoulli variables such that \eqref{mc} is satisfied. Fix $ \alpha>0,$ $\beta>0.$ Suppose that for each $\ve>0,$ there exists  $n_2\ge 1$ such that
  \begin{align}\label{al}
  &1-\ve\le \frac{\beta i}{\rho(0,i)}\le 1+\ve,\quad 1-\ve\le \frac{\beta j^{{1-\alpha}}(j^{{\alpha}}-i^{{\alpha}})}{\rho(i,j)}\le 1+\ve,
\end{align}
for $j-i\ge n_2,\ i\ge n_2.$ Then
$$\frac{{\alpha}\beta}{\log n}\sum_{j=1}^n\eta_j\overset{d}\to\xi \text{ as }n\rto,$$ where $\xi\sim \mathrm{Gamma}({\alpha},1).$
\end{theorem}
\begin{remark}
  Assume that \eqref{al} holds. If $\alpha= 1,$ then \eqref{sigma01} is naturally satisfied with $D(n)=n\beta$ for $n\ge 1.$ However, if $\alpha\ne 1,$ no function $D(\cdot)$ can satisfy \eqref{sigma01}. Specifically, for all $j\ge i\ge 1,$
\begin{align}\label{ijlu}
  (\alpha\wedge 1)(j-i)\le  j^{1-\alpha}(j^\alpha-i^\alpha)\leq(1\vee \alpha)(j-i).
\end{align}  But for any fixed $k\ge 1,$  the limit
$$\lim_{
  j=ki,\ i\rto
  }\frac{j^{{1-\alpha}}(j^{{\alpha}}-i^{{\alpha}})}{j-i}=\frac{k^{1-\alpha}(k^\alpha-1)}{k-1},$$
depends on $k$ when $\alpha\ne 1.$
\end{remark}

We outline here the key ideas for proving Theorems \ref{thbb} and \ref{tha}. To determine the limit distribution of $\sum_{i=1}^n\eta_i,$ we employ the method of moments. A primary step is to show that
for $k\ge1,$ when properly scaled, the $k$-th moment $E\z(\sum_{i=1}^n\eta_i\y)^k$ converges to a finite limit as $n\rto.$ Through detailed computations, we derive the expansion (see \eqref{exs} below) \begin{align*}
  E\Bigg(\sum_{j=1}^n\eta_j\Bigg)^k=\sum_{m=1}^k\sum_{
  \substack{l_1+\dots+l_m=k,\\
  l_s\ge1,\ s=1,...,m
  }}\frac{k!}{l_1!\cdots l_m!}\sum_{1\le j_1<...<j_m\le n}\frac{1}{\rho(0,j_1)\rho(j_1,j_2)\cdots \rho(j_{m-1},j_m)}.
  \end{align*}
    When $\rho(i,j)$ asymptotically behaves like $D(j-i)$ for some function $D(\cdot)$ and sufficiently large  $i$ and $j-i,$ the analysis reduces to studying multiple sums of the form
  \begin{align}\label{cms}
  \sum_{
    1\le j_1<...<j_m\le n
  }\frac{1}{D(j_1,s)D(j_2-j_1,s)\cdots D(j_m-j_{m-1},s)},\quad n\ge m\ge 1, \ s\in \mathbb C.
\end{align}
   In Section \ref{su}, we rigorously characterize the asymptotics of these sums. Our results yield some neat formulae
  which may be of independent interest in analysis. For cases where $\rho(i,j)$ satisfies \eqref{al}, we leverage a result from \cite{w23}, which addresses the limit behavior of related multiple sums, see Proposition \ref{thg} below.

\subsection{Applications}

%

  We employ Theorems \ref{thbb} and \ref{tha} to investigate several stochastic processes. Specific applications include determining the number of weak cutspheres for a $d (\ge 3)$-dimensional standard Brownian motion, counting the weak cutpoints of geometric Brownian motion, and analyzing the frequency with which the population size hits a predefined level in both Galton-Watson processes and branching processes with immigration in varying environments. We next describe the problems in detail.
\subsubsection{Weak cutspheres of $d(\ge3)$-dimensional Brownian motions}

Fix $d\ge 3$ and let $\{B_t\}_{t\ge0}$ be a standard $d$-dimensional Brownian motion starting at the origin.
For $x\in \mathbb R^d$ and $r>0,$ denote by
$$\mathcal B(x,r):=\{y\in \mathbb R^d: |y-x|\le r\}\text{ and } \mathcal S(x,r):=\{y\in \mathbb R^d: |y-x|=r\}$$  the ball and sphere of radius $r$ centered at $x,$ respectively.

\begin{definition}
  Fix $a>0.$ Suppose $\{B_t\}_{t\ge 0}$ is a standard $d$-dimensional Brownian motion. For $r>0,$ if the motion never returns to $\mathcal B(0,r)$ after its first visit to $\mathcal S(0,r+a),$ then we call $\mathcal S(0,r)$ a weak $a$-cutsphere of $\{B_t\}_{t\ge0}.$   For $b>a>0,$ set
\begin{align}\label{decnab}
  C(a,b):=\z\{k\in \mathbb N: \mathcal S(0,kb) \text{ is a weak } a\text{-cutsphere of } \{B_{t}\}_{t\ge 0} \y\}.
\end{align}
\end{definition}
\begin{corollary}\label{c3} Fix $d\ge 3.$  Let $\{B_t\}_{t\ge 0}$ be a standard $d$-dimensional Brownian motion. For $b>a>0,$ let $C(a,b)$ be as in \eqref{decnab}. Then
\begin{align}
  \frac{|C(a,b)\cap [1,n]|}{a/b\log n}\xrightarrow{d}\xi,
\end{align}
  as $n\to\infty,$ where $\xi\sim \mathrm{Gamma}(d-2,1).$
\end{corollary}
\begin{remark}For $1$-dimensional or $2$-dimensional Brownian motions, $C(a,b)=\emptyset$ since they are all recurrent. Therefore we need only  consider here the $d$-dimensional Brownian motions with dimension $d\ge 3.$
\end{remark}


\subsubsection{Weak cutpoints of geometric Brownian motion}
Consider the stochastic differential equation
\begin{align}
  dX_t=\mu X_tdt+\sigma X_tdB_t, \label{sde}
\end{align}
where $\mu>0,$ $\sigma\ge 0$ are constants, $\{B_t\}_{t\ge0}$ is a standard 1-dimensional Brownian motion, and $X_0=x_0>0.$ It is known that the solution to  \eqref{sde} is given by a geometric Brownian motion
\begin{align}
  X_t=X_0\exp\z\{\z(\mu-\frac{\sigma^2}{2}\y)t+\sigma B_t\y\},\ t\ge 0.\label{geob}
\end{align}
We refer to Klebaner \cite{fck} for the basics of geometric Brownian motion. Analogous to the concept of weak cutspheres for $d$-dimensional Brownian motion, we define the weak cutpoints for geometric Brownian motion as follows.
\begin{definition}
  Fix $a>0.$ Let $\{X_t\}_{t\ge 0}$ be the geometric Brownian motion in \eqref{geob}. For $x>0,$
     if $\{X_t\}_{t\ge 0}$ never returns to $(0,x]$ after its first visit to $x+a,$ then we call $x$ a weak $a$-cutpoint of $\{X_t\}_{t\ge0}.$   For $b>a>0,$ set
\begin{align}\label{decnabg}
  \hat C(a,b):=\z\{k\in \mathbb N: kb \text{ is a weak } a\text{-cutpoint of } \{X_{t}\}_{t\ge 0} \y\}.
\end{align}
\end{definition}
\begin{corollary}\label{c4}Fix $b>a>0.$ Suppose $2\mu>\sigma^2$  and let $\{X_t\}_{t\ge 0}$ be the  geometric Brownian motion given in \eqref{geob} with $x_0\in (0,b)$. Let $\hat C(a,b)$ be as in \eqref{decnabg}. Then
\begin{align}
  \frac{|\hat C(a,b)\cap [1,n]|}{a/b\log n}\xrightarrow{d}\xi,
\end{align}
  as $n\to\infty,$ where $\xi\sim \mathrm{Gamma}\z(\frac{2\mu}{\sigma^2}-1,1\y).$
\end{corollary}
\begin{remark}It is easy to see that  $X_t\to 0$ a.s. as $t\to\infty$ if $2\mu<\sigma^2,$ while $\{X_t\}_{t\ge 0}$ is recurrent if $2\mu=\sigma^2.$ Therefore, if $2\mu\le \sigma^2,$ we always have $\hat C(a,b)=\emptyset.$
\end{remark}

\subsubsection{Times the population size of a branching process reaches a given threshold}
Now, we consider a branching process with geometric offspring distribution. Set  $$f(s)=\frac{1}{2-s},\quad s\in [0,1].$$
Let $\{Y_n\}_{n\ge 0}$ be a Markov chain satisfying $P(Y_0=1)=1$ and
\begin{align*}
  E\z(s^{Y_n}\,\middle |\,Y_{n-1},\ldots,Y_{0}\y)=f(s)^{Y_{n-1}},\quad \ s\in [0,1],\ n\ge 1.
\end{align*}
This defines a critical Galton-Watson branching process with geometric offspring distribution. 
For  $n\ge1,$ define \begin{align*}
 N_n:=\#\{1\le t\le n:Y_t=1\}
\end{align*}
  as the number of generations where the population size is exactly $1.$ The following corollary characterizes the limit distribution of $N_n$.
\begin{corollary}\label{thy}
  We have
  $N_n\xrightarrow{a.s.}\xi$
 as $n\rto,$ where $\xi\sim\mathrm{Geo}\z(\frac{6}{\pi^2}\y)$.
\end{corollary}

Next, we consider a branching process in varying environments (BPVE) with immigration. It is well-known that compared to Galton-Watson processes, BPVEs exhibit new phenomena due to the influence of time-dependent environments. For foundational studies on this topic, we refer to  \cite{fuj,ker,kv,lin,ms,s24,wb24}. To formalize the model,
let $0<p_i<1$ for $i\ge 1$ be a sequence of numbers. For $i\ge 1$, define the probability generating function
\begin{align}\label{fis}
   \tilde f_i(s)=\frac{p_i}{1-(1-p_i)s},\quad s\in [0,1],
 \end{align}
which governs the offspring distribution of an individual from the $(i-1)$-th generation.
Let $\{Z_n\}_{n\ge 0}$ be a Markov chain satisfying $P(Z_0=0)=1$ and
 \begin{align}
   E\z(s^{Z_{n}}\,\middle |\,Z_0,\ldots,Z_{n-1}\y)=\tilde f_n(s)^{Z_{n-1}+1},\quad n\ge 1. \label{mb}
 \end{align}
 This defines a BPVE with exactly one immigrant per generation.
 It follows from \eqref{mb} that
      $$m_i:=E\z(Z_i\,\middle |\,Z_{i-1}=0\y)=\tilde f_i'(1)=\frac{1-p_i}{p_i},\quad i\ge 1.$$
  We focus on the near-critical regime, for which $m_i\to 1$ (equivalently $p_i\to 1/2$) as $i\to\infty.$
 Specifically, set
 \begin{align}\label{dpa}
   p_i=
   \frac{1}{2}-\frac{r_i}{4},\quad i\ge 1,
 \end{align} where $r_i\in [0,1]$ for all $\ i\ge1.$
   For $n\ge 1,$ define
$$I_n=\#\{1\le t\le n: Z_t=0\}$$ as the number of regenerating times within $[1,n].$ A time $t$ is termed a regenerating time if
$Z_t=0.$
On such events, all individuals (including the immigrant) from generation
 $t-1$ leave no descendants, leading to local extinction.
 However, an immigrant arrives almost surely in generation $t,$ thereby regenerating the process. Thus, $I_n$ quantifies the frequency of such regenerative events.

 \begin{corollary} \label{thz}
   {\rm(i)} Suppose $p_i,\ i\ge 1$ are defined by \eqref{dpa} and $\lim_{i\to\infty}\sum_{n\ge i}r_n=0$. Then
 $\frac{I_n}{\log n}\xrightarrow{d}\xi$ as $n\rto,$
 where $\xi\sim\mathrm{Exp}(1)$.
 {\rm (ii)} Fix $B\in [0,1).$ If $p_i=\frac{1}{2}+\frac{B}{4i}$ for $i\ge1,$ then
 $\frac{I_n}{\log n}\xrightarrow{d}\xi$  as $n\to\infty,$
 where $\xi\sim\mathrm{Gamma}(1-B,1)$.
  \end{corollary}
\begin{remark}\label{rmz} 
In \cite{wb24}, the first author studies the case $p_i=\frac{1}{2}+\frac{B}{4i},\; i\ge1 .$ Fix a positive integer $a,$ and define  $C(a)=\#\{t\ge 0: Z_t=a\}.$ It is shown that the set $C(a)$  is almost surely infinite if $0\le B<1,$ and finite if $B\ge 1.$ Specifically, when $B=0,$ it is further  proved in \cite{wb24} that
$\#\{C(a)\cap[1,n]\}/\log n\xrightarrow{d}\xi$ as $n\rto,$ where $\xi\sim \mathrm{Exp}(1).$
   \end{remark}

\noindent\textbf{Outline of the paper.} The remainder of the paper is organized as follows. As mentioned above, Section \ref{su} investigates the asymptotic behavior of the multiple sums defined in \eqref{cms}. Building on these results, Section \ref{sb} establishes the proofs of Theorems \ref{thbb} and \ref{tha}. The proofs of Corollaries \ref{c3}, \ref{c4}, \ref{thy} and \ref{thz}, which are derived from Theorems \ref{thbb} and \ref{tha}, are detailed in Section \ref{sbp}. Finally, concluding remarks are presented in Section \ref{sc}.



\section{Asymptotic analysis of multiple sums}\label{su}
%

A key step in proving Theorem \ref{thbb} involves studying the asymptotic behavior of the multiple sums defined in \eqref{cms}. The following theorem characterizes the limiting behavior of these sums across distinct regimes.

\begin{theorem}\label{prpd}
 {\rm(i)} Fix $n_0\ge 1.$ For $s\in \mathbb C,$ let $D(n,s),\ n\ge 1,$ be complex functions such that $\sum_{n=1}^\infty\frac{1}{|D(n,s)|}<\infty$ and denote $\lambda(s):=\sum_{n=n_0}^\infty \frac{1}{D(n,s)}.$
     For $1\le m\le n$  and $s\in \mathbb C,$ set \begin{align*}
  \Phi(n,m,s)=\sum_{\substack{
    1\le j_1<...<j_m\le n\\
    j_i-j_{i-1}\ge n_0,\ i=1,...,m
  }}\frac{1}{D(j_1,s)D(j_2-j_1,s)\cdots D(j_m-j_{m-1},s)},
\end{align*}
where and throughout, we assume $j_0=0.$
Then, we have \begin{align}\label{lgnm}
  \lim_{n\rto}\Phi(n,m,s)=\lambda(s)^m.
\end{align}
{\rm(ii)} Fix $n_0\ge 1.$ Let $D(n),\ n\ge1$ be a sequence of numbers such that $\inf_{n\ge1}D(n)>0$  and $\sum_{n=1}^\infty \frac{1}{D(n)}=\infty.$ For $n\ge 1,$ set $S(n):=\sum_{i=1}^n \frac{1}{D(i)}.$ For $m\ge 1,$ denote \begin{align*}
  \Phi(n,m)=\sum_{\substack{
    1\le j_1<...<j_m\le n\\
    j_i-j_{i-1}\ge n_0,\ 1\le i\le m}}\frac{1}{D(j_1)D(j_2-j_1)\cdots D(j_m-j_{m-1})}.
\end{align*}
  If $\{S(n)\}_{n\ge1}$ is regularly varying with index $\tau\in [0,1],$
then we have
\begin{align}\label{lgnmd}
  \lim_{n\rto}\frac{\Phi(n,m)}{(S(n))^m}=\lambda_\tau^{-(m-1)},
\end{align}
where $\lambda_0=1$ and $\lambda_\tau=\frac{\Gamma(1+2\tau)}{\tau\Gamma(1+\tau)\Gamma(\tau)}$ for $\tau\in (0,1].$
\end{theorem}
Based on Theorem \ref{prpd}, we derive some formulae that may have potential applications in analysis. To formalize these results, we first introduce necessary notations. Let $\log_0i=i,$ and for  $n\ge1,$ set
$\log_n i=\log\z(\log_{n-1}i\y).$
For an integer $m\ge 0$ and  a complex number $s,$ set
\begin{align}\label{dom}
  \mathscr O_m=\min\z\{n\in \mathbb N: \log_m n>0\y\},\ \lambda(m,s,i)=\z(\log_mi\y)^s\prod_{j=0}^{m-1} \log_{j}i,
\end{align}
where and in what follows, we assume that an empty product equals 1 and an empty sum equals $0.$
Now fix integers $m\ge0$ and $k\ge1.$ For $s=\sigma+t\mathrm i\in \mathbb C$ and $n_0\ge \mathscr O_m,$ define
\begin{align}\label{defg}
    &U_n(k,m,n_0, s)=\sum_{\substack{
    1\le j_1<...<j_k\le n\\
    j_i-j_{i-1}\ge n_0,\ i=1,...,k
  }}\frac{1}{\lambda(m,s,j_1)\lambda(m,s,j_2-j_1)\cdots \lambda(m,s,j_{k}-j_{k-1})}.
\end{align}

\begin{corollary}\label{rzr}
  Fix integers $m\ge0,$ $k\ge1$ and let $\mathscr O_m$ be as in \eqref{dom}. For $n_0\ge \mathscr O_m$ and $s=\sigma+t\mathrm i\in \mathbb C,$ let $U_n(k,m,n_0, s)$ be defined as in \eqref{defg}.
  \begin{itemize}
    \item[{\rm(i)}] If $\sigma>1,$ then
  \begin{align}
    \lim_{n\rto} U_n(k,m,n_0, s)=\zeta(m,s)^k, \label{rzf}
  \end{align}
  where $\zeta(m,s)=\sum_{i=n_0}^\infty \frac{1}{\lambda(m,s,i)}$ is absolutely convergent.
   \item[{\rm(ii)}] If $\sigma=1,$ then
  \begin{align*}
    \lim_{n\rto} \frac{U_n(k,m,n_0,\sigma)}{\z(\log_{m+1}n\y)^k}=1.
  \end{align*}
  \item[{\rm(iii)}] If $0\le \sigma<1$ and $m\ge1,$ then
  \begin{align*}
    \lim_{n\rto} \frac{(1-\sigma)^kU_n(k,m,n_0,\sigma)}{(\log_{m}n)^k}=1.
  \end{align*}
  \item[{\rm(iv)}] If $0\le\sigma<1$ and $m=0,$ then
   \begin{align}\label{sl1}
   &\lim_{n\rto}\frac{U_n(k,m,n_0,\sigma)}{n^{k(1-\sigma)}}=\frac{1}{1-\sigma}\z(\frac{\Gamma(2-\sigma)\Gamma(1-\sigma)}{\Gamma(3-2\sigma)}\y)^{k-1}. \end{align}
  \end{itemize}

\end{corollary}

\begin{remark} (a) We shed here some light on the case $m=0.$ Here $\lambda(0,s,i)=i^s$ and $\mathscr O_0=1.$ Consider $$\zeta(s)=\sum_{n=1}^\infty \frac{1}{n^s},\quad s\in \mathbb C,$$ which is known as the Riemann zeta function and is  absolutely convergent if $\mathrm{Re}(s)>1.$ Setting $n_0=\mathscr O_0=1,$
Part (i) of Corollary \ref{rzr} yields that
\begin{align*}
    &\lim_{n\to\infty}\sum_{\begin{subarray}{c}
    1\le j_1<...<j_k\le n
  \end{subarray}}\frac{1}{j_1^s(j_2-j_1)^s(j_3-j_2)^{s}\cdots (j_{k}-j_{k-1})^s}=\zeta(s)^k,\quad k\ge1,\  \mathrm{Re}(s)>1.
\end{align*}
  This provides an alternative representation of the Riemann zeta function.
 We remark that the function $\zeta(s),$ initially defined for $s\in \mathbb C$ such that $\mathrm{Re}(s)>1,$ admits an analytic continuation to the entire complex plane which is holomorphic except for a simple pole at $s=1$ with residue $1.$ For further details, see \cite{Ed,kavo}.

 Furthermore,  from \eqref{sl1},  for $0<\sigma<1$ and $n_0\ge 1,$ we have
  \begin{align*}
 \lim_{n\rto}\frac{1}{n^{k(1-\sigma)}}&\sum_{\substack{
   1\le j_1<...<j_k\le n\\
   j_s-j_{s-1}\ge n_0
 }}\frac{1}{j_1^\sigma(j_2-j_1)^\sigma\cdots (j_{m}-j_{m-1})^\sigma}=\frac{1}{1-\sigma}\z(\frac{\Gamma(2-\sigma)\Gamma(1-\sigma)}{\Gamma(3-2\sigma)}\y)^{k-1}.
  \end{align*}

 (b) Notice that in Parts (ii)--(iv) of Corollary \eqref{rzr}, we consider $U_n(k,m,\sigma)$ with $\sigma\in \mathbb R,$ rather than the complex variable $s$. This restriction arises because it makes no sense to consider the asymptotics of  $U_n(k,m,s)$ for $s\in \mathbb C$ with $0<\text{Re}(s)\le 1.$ For example, let $s=\sigma+\mathrm{i}t\in \mathbb C.$
  Then by definition  we have
$$U_n(1,m,s)=\sum_{j=\mathscr O_m}^n \frac{\cos \z(t\log_{m+1} j\y)}{\lambda(m,\sigma, j)}-\mathrm{i}\frac{\sin \z(t\log_{m+1} j\y)}{\lambda(m,\sigma, j)}.$$
 For $t\ne 0$ and $0<\sigma\le1,$
neither the real part $\sum_{j=1}^{n}\cos \z(t\log_{m+1} j\y)/\lambda(m,\sigma, j)$ nor the imaginary part  $\sum_{j=1}^n\sin \z(t\log_{m+1} j\y)/j^\sigma$ converges as $n\rto$. Thus, for $\sigma\in (0,1],$ the asymptotics of $U_n(k,m,s)$ are meaningful only when $t=0.$

\end{remark}
Next, we quote a result from \cite{w23}, which will be used to  prove Theorem \ref{tha}.
    \begin{proposition}[\cite{w23}]\label{thg}
   Fix a number ${\alpha}>0$ and an integer $k\ge 1$. Suppose $b_s,$ $s=1,..,k$ are  positive integers and let $j_0=0.$ Then,
\begin{align}
  \lim_{n\rto}\frac{1}{(\log n)^k}{\sum_{\substack{
    1\le j_1<...<j_k\le n\\
    j_s-j_{s-1}\ge b_s,\, s=1,...,k}
}\frac{1}{j_1j_2^{1-\alpha}(j_2^{{\alpha}}-j_1^{{\alpha}})\cdots j_k^{1-\alpha}(j_k^{{\alpha}}-j_{k-1}^{{\alpha}})}}=\frac{\prod_{j=0}^{k-1}(j+\alpha)}{k!\alpha^k}.\label{aln}
\end{align}
 \end{proposition}
We divide the rest of this section into three subsections. Parts (i) and (ii) of Theorem \ref{prpd} are proved in Subsections \ref{susec21} and \ref{susec22} respectively. Corollary \ref{rzr} is proved in Subsection \ref{subsec23}.

\subsection{Proof of Part (i) of Theorem \ref{prpd}}\label{susec21}

\proof We prove Part (i) of Theorem \ref{prpd} by induction.  Since  $\sum_{j_1=n_0}^\infty \frac{1}{D(j_1,s)}=\lambda(s),$ \eqref{lgnm} holds trivially for $m=1.$ Fix $k\ge 2$ and suppose \eqref{lgnm} holds for $m=k-1.$ We now prove it also holds for $m=k.$ It is easy to see that for $1\leq k\leq n$,
\begin{align}\label{gkk1}
  \Phi&(n,k,s)=\sum_{\substack{
    1\le j_1<...<j_k\le n\\
    j_i-j_{i-1}\ge n_0,\, i=1,...,k}}\frac{1}{D(j_1,s)D(j_2-j_1,s)\cdots D(j_k-j_{k-1},s)}\\
  &=\sum_{j_1=n_0}^{n-(k-1)n_0}\frac{1}{D(j_1,s)}\sum_{\substack{
     j_1<j_2<...<j_k\le n\\
    j_i-j_{i-1}\ge n_0,\, i=2,...,k}}\frac{1}{D(j_2-j_1,s)D(j_3-j_2,s)\cdots D(j_k-j_{k-1},s)}\no\\
 & =\sum_{j_1=n_0}^{n-(k-1)n_0}\frac{1}{D(j_1,s)}\sum_{\substack{
   1\le j_1'<...<j_{k-1}'\le n-j_1\\
   j_i-j_{i-1}\ge n_0, \,i=1,...,k-1}}\frac{1}{D(j_1',s)D(j_2'-j_1',s)\cdots D(j_{k-1}'-j_{k-2}',s)}\no\\
 &=\sum_{j_1=n_0}^{n-(k-1)n_0}\frac{1}{D(j_1,s)}\Phi(n-j_1,k-1,s)\no\\
 &=\left[\sum_{j_1=n_0}^{\lfloor n/2\rfloor}+\sum_{j_1=\lfloor n/2\rfloor+1}^{n-(k-1)n_0}\right]\frac{1}{D(j_1,s)}\Phi(n-j_1,k-1,s)\no\\
 &=: {\rm I_n(s)+II_n(s)}.\no
\end{align}

Consider first the term $\mathrm I_n(s)$ in \eqref{gkk1}.
Observe that
\begin{align}\label{ione}
  \mathrm I_n(s)&=\sum_{j_1=n_0}^{\lfloor n/2\rfloor}\frac{1}{D(j_1,s)}\Phi(n-j_1,k-1,s)\\
  &=\sum_{j_1=n_0}^{\lfloor n/2\rfloor}\frac{1}{D(j_1,s)}\lambda(s)^{k-1}+\sum_{j_1=n_0}^{\lfloor n/2\rfloor}\frac{1}{D(j_1,s)}\z(\Phi(n-j_1,k-1,s)-\lambda(s)^{k-1}\y).\no
\end{align}
Fix $\ve>0.$ Since
$\lim_{n\rto}\Phi(n,k-1,s)=\lambda(s)^{k-1},$   there exists $N>0$ such that
\begin{align*}
  |\Phi(n,k-1,s)-\lambda(s)^{k-1}|<\ve, \text{ for all }n\ge N.
\end{align*}
Thus, for $n\ge 2N,$ we have
\begin{align}
\z| \sum_{j_1=n_0}^{\lfloor n/2\rfloor}\frac{1}{D(j_1,s)}\z(\Phi(n-j_1,k-1,s)-\lambda(s)^{k-1}\y)\y|\le \ve \sum_{j_1=n_0}^{\lfloor n/2\rfloor}\frac{1}{|D(j_1,s)|}\le \ve \sum_{j_1=n_0}^{\infty}\frac{1}{|D(j_1,s)|}.\no
\end{align}
Therefore, we deduce that
\begin{align*}
   \lim_{n\rto}\sum_{j_1=n_0}^{\lfloor n/2\rfloor}\frac{1}{D(j_1,s)}\z(\Phi(n-j_1,k-1,s)-\lambda(s)^{k-1}\y)=0.\no
\end{align*}
Consequently, letting $n\rto$ in \eqref{ione}, we get
\begin{align}\label{lp1}
  \lim_{n\rto}\mathrm I_n(s)=\lambda(s)^{k}.
\end{align}

Consider next the term $\mathrm{II}_n(s)$ in \eqref{gkk1}.  Since
$\lim_{n\rto}\Phi(n,k-1,s)=\lambda(s)^{k-1},$ there exists $M>0$ such that
$|\Phi(n,k-1,s)|<M$ for all $n\ge1.$ Consequently,
\begin{align*}
  |\mathrm{II}_n(s)|&\le \sum_{j_1=\lfloor n/2\rfloor+1}^{n-(k-1)n_0}\frac{1}{|D(j_1,s)|}|\Phi(n-j_1,k-1,s)|\le M \sum_{j_1=\lfloor n/2\rfloor+1}^{n-(k-1)n_0}\frac{1}{|D(j_1,s)|}.
\end{align*}
Recalling that $\sum_{n=1}^\infty\frac{1}{|D(n,s)|}<\infty,$  thus we have
\begin{align}
  \lim_{n\rto}|\mathrm{II}_n(s)|\le M\lim_{n\rto}\sum_{j_1=\lfloor n/2\rfloor+1}^{n-(k-1)n_0}\frac{1}{|D(j_1,s)|}=0. \label{lp2}
\end{align}
Plugging \eqref{lp1} and \eqref{lp2} into \eqref{gkk1}, we conclude that
 $\lim_{n\rto}\Phi(n,k, s)={\lambda(s)}^{k}.$ Thus, by induction,   \eqref{lgnm} holds for all $k\ge1.$ This completes the proof of Part (i) of Theorem  \ref{prpd}.  \qed

\subsection{Proof of Part (ii) of Theorem \ref{prpd}} \label{susec22}
\proof By assumption, we have
$\inf_{n\ge1}D(n)>0$ and $\sum_{n=1}^\infty \frac{1}{D(n)}=\infty.$ Consequently,
\begin{align}\label{uubd}
  \lim_{n\rto}S(n)=\infty \text{ and } 0<\sup_{n\ge1}{1}/{D(n)}<C,
\end{align}  for some $C>0.$ Throughout the remainder of the proof, we will implicitly utilize the facts stated in \eqref{uubd} without further explicit reference.

 By virtue of \eqref{uubd}, it follows that
 $$\lim_{n\rto}\frac{\Phi(n,1)}{S(n)}=\lim_{n\rto}\sum_{j_1=n_0}^n\frac{1}{D(j_1)}\bigg/\sum_{i=1}^n\frac{1}{D(i)}=1,$$
 which confirms that \eqref{lgnmd} holds for $m=1.$
 For $k\ge1,$ a recursive argument analogous to \eqref{gkk1} yields that
 \begin{align}\label{pki}
  \Phi(n,k+1)&=\sum_{\substack{
    1\le j_1<...<j_{k+1}\le n\\
    j_i-j_{i-1}\ge n_0,\, 1\le i\le k+1}}\frac{1}{D(j_1)D(j_2-j_1)\cdots D(j_{k+1}-j_{k})}\\
  &=\sum_{j_1=n_0}^{n-kn_0}\frac{1}{D(j_1)}\Phi(n-j_1,k)\no\\
  &=\left[\sum_{j_1=n_0}^{\lfloor n/2\rfloor}+\sum_{j_1=\lfloor n/2\rfloor+1}^{n-kn_0}\right]\frac{1}{D(j_1)}\Phi(n-j_1,k).\no
\end{align}
Fix $k\ge 1,$ and assume \eqref{lgnmd} holds for $m=k.$ We now prove it holds  for $m=k+1.$ The proof is divided into two cases: $\tau=0$ and $\tau\in (0,1].$

First, assume that $\tau=0.$  Since $\Phi(n,k+1)$ is increasing in $n,$ we obtain
 \begin{align}
   \sum_{j_1=n_0}^{\lfloor n/2\rfloor}\frac{1}{D(j_1)}\Phi(n-\lfloor n/2\rfloor,k)\le \Phi(n,k+1)\le \sum_{j_1=n_0}^{n}\frac{1}{D(j_1)}\Phi(n,k).\label{lupk}
\end{align}
 Given the inductive hypothesis and the regular variation of
$\{S(n)\}_{n\ge1}$ with index $\tau,$
  we have $$\lim_{n\rto}\frac{\sum_{j_1=n_0}^{\lfloor n/2\rfloor}\frac{1}{D(j_1)}}{S(n)}=\lim_{n\rto}\frac{S(\lfloor n/2\rfloor)}{S(n)}=1 $$
and
$$\lim_{n\rto}\frac{\Phi(n-\lfloor n/2 \rfloor,k)}{(S(n))^k}=\lim_{n\rto}\frac{\Phi(n,k)}{(S(n))^k}=1.$$
Dividing \eqref{lupk} by $(S(n))^{k+1}$ and taking $n\to\infty$, we obtain
$$\lim_{n\rto}\frac{\Phi(n,k+1)}{(S(n))^{k+1}}=1.$$
 Thus by induction, \eqref{lgnmd} holds for $m=k+1.$


Next, assume $\tau \in (0, 1]$. Fix an integer $l \ge 2$. From \eqref{pki} we have
\begin{align*}
\Phi(n,k+1) &= \sum_{j_1=n_0}^{n-kn_0}\frac{1}{D(j_1)}\Phi(n-j_1,k) \nonumber \\
&= \sum_{s=1}^{l-1}\sum_{i=(s-1)u+n_0+1}^{su+n_0}\frac{1}{D(i)}\Phi(n-i,k) \nonumber \\
&\quad + \sum_{i=(l-1)u+n_0+1}^{n-kn_0}\frac{1}{D(i)}\Phi(n-i,k) + \frac{\Phi(n-n_0,k)}{D(n_0)},
\end{align*}
where $u = \lfloor(n - (k+1)n_0)/l \rfloor$ is temporarily defined. Since $\Phi(n,m)$ is monotonically increasing in $n$, we obtain
\begin{align} \label{pkul}
  \Phi&(n,k+1)\le \sum_{s=1}^{l-1}\sum_{i=(s-1)u+n_0+1}^{su+n_0}\frac{1}{D(i)}\Phi(n-(s-1)u-n_0-1,k)\\
&\quad\quad\quad\quad\quad\quad+\sum_{i=(l-1)u+n_0+1}^{n-kn_0}\frac{1}{D(i)}\Phi(n-(l-1)u-n_0-1,k)+\frac{\Phi(n-n_0,k)}{D(n_0)}\no\\
    &=  \frac{\Phi(n-n_0,k)}{D(n_0)}+\sum_{s=1}^{l-1}\z[S(su+n_0)-S((s-1)u+n_0)\y]\Phi(n-(s-1)u-n_0-1,k)\no\\
  &\quad\quad+ \z[S(n-kn_0)-S((l-1)u+n_0)\y]\Phi(n-(l-1)u-n_0-1,k). \no
 \end{align}
Recall that \eqref{lgnmd} is assumed to be true for $m=k.$ Dividing \eqref{pkul} by $S(n)^{k+1}$ and taking the upper limit, we obtain
\begin{align}\label{uul}
  \varlimsup_{n\rto}\frac{\Phi(n,k+1)}{\z(S(n)\y)^{k+1}}&\le \lambda_\tau^{-(k-1)}\sum_{s=1}^{l-1}\z[\z({s}/{l}\y)^{\tau}-\z(({s-1})/{l}\y)^{\tau}\y]\z(({l-s+1})/{l}\y)^{\tau}\\  &\quad\quad\quad\quad+ \lambda_\tau^{-(k-1)}\z[\z({l}/{l}\y)^{\tau}-\z(({l-1})/{l}\y)^{\tau}\y]\z({1}/{l}\y)^{\tau}\no\\
  &= \lambda_\tau^{-(k-1)}\sum_{s=1}^{l}\z[\z({s}/{l}\y)^{\tau}-\z(({s-1})/{l}\y)^{\tau}\y]\z[{(l-s+1)}/{l}\y]^{\tau}.\no
\end{align}
For the lower limit, similar arguments yield
\begin{align}
  \Phi(n,k+1)&= \sum_{j_1=n_0}^{n-kn_0}\frac{1}{D(j_1)}\Phi(n-j_1,k)
  \ge\sum_{s=1}^{l-1}\sum_{i=(s-1)u+n_0+1}^{su+n_0}\frac{1}{D(i)}\Phi(n-su-n_0,k)\no\\
    &= \sum_{s=1}^{l-1}\z[S(su+n_0)-S((s-1)u+n_0)\y]\Phi(n-su-n_0,k).\no
\end{align}
Thus, using \eqref{lgnmd} for $m=k,$ we deduce that
\begin{align}\label{lul}
 \varliminf_{n\rto}\frac{\Phi(n,k+1)}{\z(S(n)\y)^{k+1}}&\ge \lambda_\tau^{-(k-1)}\sum_{s=1}^{l-1}\z[\z({s}/{l}\y)^{\tau}-\z({(s-1)}/{l}\y)^{\tau}\y]\z(({l-s})/{l}\y)^{\tau}.  \end{align}
Combining \eqref{uul} and \eqref{lul} gives
\begin{align}\label{uslu}
 \lambda_\tau^{-(k-1)}&\sum_{s=1}^{l-1}[\z({s}/{l}\y)^{\tau} -\z({(s-1)}/{l}\y)^{\tau}] \z(({l-s})/{l}\y)^{\tau}\le \varliminf_{n\rto}\frac{\Phi(n,k+1)}{\z(S(n)\y)^{k+1}}\\
  &
  \le \varlimsup_{n\rto}\frac{\Phi(n,k+1)}{\z(S(n)\y)^{k+1}}\le
   \lambda_\tau^{-(k-1)}\sum_{s=1}^{l}\z[\z({s}/{l}\y)^{\tau}-\z(({s-1})/{l}\y)^{\tau}\y]\z[({l-s+1})/{l}\y]^{\tau}.\no
\end{align}
Next, we demonstrate that
\begin{align}\label{ubl}
  \lim_{l\rto}&\frac{1}{l^{2\tau}}\sum_{s=1}^{l}\z(s^{\tau}-\z(s-1\y)^{\tau}\y)\z(l-s+1\y)^{\tau}= \lambda_\tau^{-1},\\
\label{lbl}
  \lim_{l\rto}&\frac{1}{l^{2\tau}}\sum_{s=1}^{l-1}\z(s^{\tau}-\z(s-1\y)^{\tau}\y)\z(l-s\y)^{\tau}= \lambda_\tau^{-1}.
\end{align}
 For $\tau=1,$ observe that $\lambda_1=\frac{\Gamma(3)}{\Gamma(2)\Gamma(1)}=2.$ Direct computation confirms that both \eqref{ubl} and \eqref{lbl} hold in this case.

Suppose next $\tau\in (0,1),$ and
 write
$$g(s)=\z(s^\tau-\z(s-1\y)^\tau\y)\z(l-s+1\y)^\tau,\quad s\in[1,l].$$
Since both $s^\tau-\z(s-1\y)^\tau$ and $(l-s+1)^\tau$ are decreasing in $s$, $g(s)$ is also decreasing on $[1,l]$. Fix $\varepsilon>0.$ By the mean value theorem, there exists $\eta\in(s-1,s)$ such that for $s\geq1$,
\begin{align}
 \tau s^{\tau-1}\le s^{\tau}-(s-1)^{\tau}=\tau\eta^{\tau-1}<\tau(s-1)^{\tau-1}.\no
\end{align}
Combining this with the monotonicity of $g(s)$, we infer that
\begin{align}
  \varlimsup_{l\rto}&\frac{1}{l^{2\tau}}\sum_{s=1}^{l}\z(s^{\tau}-\z(s-1\y)^{\tau}\y)\z(l-s+1\y)^{\tau}\no\\
  &=\varlimsup_{l\rto}\frac{1}{l^{2\tau}}\z[\int^{l}_{1}\z(s^{\tau}-\z(s-1\y)^{\tau}\y)\z(l-s+1\y)^{\tau}ds+l^\tau\y]\no\\
  &\le \tau\varlimsup_{l\rto}\frac{1}{l^{2\tau}}\int_{1}^l (s-1)^{\tau-1}(l-s+1)^rds\no\\
  &=\tau\varlimsup_{l\rto}\int_{1/l}^1 \z(x-\frac{1}{l}\y)^{\tau-1}\z(1-x+\frac{1}{l}\y)^{\tau}dx\no\\
    &=\tau\int_{0}^1 x^{\tau-1}\z(1-x\y)^{\tau}dx=\frac{\tau\Gamma(\tau+1)\Gamma(\tau)}{\Gamma(2\tau+1)}=\lambda_\tau^{-1}.\no
\end{align}
Similarly, we can show that
\begin{align}
  \varliminf_{l\rto}&\frac{1}{l^{2\tau}}\sum_{s=1}^{l}\z(s^{\tau}-\z(s-1\y)^{\tau}\y)\z(l-s+1\y)^{\tau}\ge \lambda_\tau^{-1}.\no
\end{align}
Thus, \eqref{ubl} holds. An analogous argument confirms \eqref{lbl}.

 Since \eqref{uslu} holds for all $l\ge 2,$ taking $l\rto$ and applying \eqref{ubl} and \eqref{lbl}, we conclude that
 \begin{align}
   \lim_{n\rto}\frac{\Phi(n,k+1)}{\z(S(n)\y)^{k+1}}= \lambda_\tau^{-k},\no
 \end{align} which establishes  \eqref{lgnmd} for $m=k+1.$
 By induction, \eqref{lgnmd} holds for all $k\ge1.$  Part (ii) of Theorem \ref{prpd} is proved. \qed


\subsection{Proof of Corollary \ref{rzr}}\label{subsec23}

\proof Fix   $m\in \mathbb N$  and  $0<\sigma\in\mathbb R.$  Let $\mathscr O_m$ be as in \eqref{dom}.  For $k\ge \mathscr O_m$ and $s=\sigma+\mathrm it\in \mathbb C,$ set $D(k,s)=\lambda(m,s,k)=\z(\log_mk\y)^s\prod_{j=0}^{m-1}\log_j k.$
For $x\ge \mathscr O_m,$ let $A(x):=\int_{\mathscr O_m}^x\frac{1}{D(u,\sigma)}du.$
A direct computation yields
\begin{align*}
A(x)=\int_{\mathscr O_m}^x\frac{1}{\lambda(m,\sigma,u)}du=
\left\{\begin{array}{ll}
  \frac{1}{1-\sigma}\z(\z(\log_m x\y)^{1-\sigma}-\z(\log_m \mathscr O_m\y)^{1-\sigma}\y), & \text{if } \sigma\ne 1,\\
  \log_{m+1}x-\log_{m+1}\mathscr O_m, &\text{if } \sigma=1.
\end{array}\right. \end{align*}

If $\sigma>1,$ the series $\sum_{n=\mathscr O_m}^\infty\frac{1}{D(n,s)}$ is absolutely convergent. Applying Part (i) of Theorem \ref{prpd} we get  \eqref{rzf}.

For $n\ge1,$ let $$D(n)=\z\{\begin{array}{ll}
D(n,\sigma), & \text{if }n\ge \mathscr O_m,\\
D(\mathscr O_m,\sigma),&\text{if }1\le n< \mathscr O_m,
\end{array}\y. $$  and set $S(n)=\sum_{i=1}^n\frac{1}{D(i)}.$

If $0\le \sigma\le 1,$ it is easy to see that $\inf_{n\ge1}D(n)\ge \z(\log_m \mathscr O_m\y)^s\prod_{j=0}^{m-1}\log_j  \mathscr O_m>0$ and
$S(n)\sim A(n)\to \infty$ as $n\rto.$ Thus for $x>0,$ we have
\begin{align*}
 \lim_{n\to\infty}\frac{S(\lfloor xn\rfloor)}{S(n)} =x^{\tau} \text{ with } \tau=\z\{\begin{array}{ll}
    0, & \text{if } m\ge1,\\
    1-\sigma,&\text{if }m=0.
  \end{array}\right.
\end{align*}
 Consequently, the sequence $\{S(n)\}_{n\ge1}$ is regularly varying  with index $\tau.$
Applying Part (ii) of Theorem \ref{prpd},  we obtain Parts (ii)--  (iv) of Corollary \ref{rzr}.
\qed


\section{Proofs of Theorems \ref{thbb} and \ref{tha}}\label{sb}
In this section, we prove Theorems \ref{thbb} and \ref{tha}  by the method of moments. The key step involves computing the moments  $E\z(\sum_{j=1}^n\eta_j\y)^k$ for $k\ge1.$

Fix $k\ge1.$  Applying the multinomial expansion theorem, we get
\begin{align}
  E\Bigg(\sum_{j=1}^n\eta_j\Bigg)^k&=\sum_{\begin{subarray}{c}
i_1+...+i_n=k,\\
i_s\ge 0,\, s=1,...,n
  \end{subarray}}\frac{k!}{i_1!i_2!\cdots i_n!}E\z(\eta_1^{i_1}\eta_{2}^{i_2}\cdots \eta_{n}^{i_n}\y)\no\\
&=\sum_{m=1}^k\sum_{1\le j_1<...<j_m\le n}\sum_{\begin{subarray}{c}
i_{j_1}+...+i_{j_m}=k\\
i_{j_s}>0,\, s=1,...,m
  \end{subarray}}\frac{k!}{i_{j_1}!i_{j_2}!\cdots i_{j_m}!}E\z(\eta_{j_1}^{i_{j_1}}\eta_{j_2}^{i_{j_2}}\cdots \eta_{j_m}^{i_{j_m}}\y),\no
  \end{align}
where we adopt the convention $0^0=1.$
  Since $\eta_i\in \{0,1\}$ for all $i\ge 1,$ taking \eqref{mc} into account, we have
\begin{align}\label{exs}
  E\Bigg(\sum_{j=1}^n\eta_j\Bigg)^k&=\sum_{m=1}^k\sum_{\begin{subarray}{c}
l_1+...+l_m=k\\
l_{s}>0,\,s=1,...,m
  \end{subarray}}\frac{k!}{l_1!l_{2}!\cdots l_{m}!}\sum_{1\le j_1<...<j_m\le n}E\Big(\eta_{j_1}\eta_{j_2}\cdots \eta_{j_m}\Big)\\
    &=\sum_{m=1}^k\sum_{\begin{subarray}{c}l_1+\dots+l_m=k,\\
  l_s\ge1,\,s=1,...,m
  \end{subarray}}\frac{k!}{l_1!\cdots l_m!}\sum_{1\le j_1<...<j_m\le n}P\z(\eta_{j_1}=1,\eta_{j_2}=1,\cdots ,\eta_{j_m}=1\y)\no\\
  &=\sum_{m=1}^k\sum_{\begin{subarray}{c}l_1+\dots+l_m=k,\\
  l_s\ge1,\,s=1,...,m
  \end{subarray}}\frac{k!}{l_1!\cdots l_m!}\Psi_n(m),\no
  \end{align}
  where and in what follows
  \begin{align}\label{depsi}
    \Psi_n(m)=\sum_{1\le j_1<...<j_m\le n}P\z(\eta_{j_1}=1,\eta_{j_2}=1,\cdots ,\eta_{j_m}=1\y),\quad n\ge m \ge 1.
  \end{align}

  Based on \eqref{exs},  we divide the remainder of this section into three subsections, which will finish the proofs of Parts (i)--(ii) of Theorem \ref{thbb} and Theorem \ref{tha}, respectively.
\subsection{Proof of Part (i) of Theorem \ref{thbb}}
  \proof Let $D(n),\ n\ge 0$ be a sequence satisfying $D(0)=1,$ $D(n)>1$ for  $n\ge1$ and \begin{align}
\zeta(D)= \sum_{n=1}^\infty \frac{1}{D(n)}<\infty.\label{zed}
  \end{align}By assumption, we have
  \begin{align}
    P(\eta_i=1)=\frac{1}{\rho(0,i)}, \quad P\z(\eta_j=1\middle|\,\eta_i=1\y)=\frac{1}{\rho(i,j)}=\frac{1}{D(j-i)},\label{dns}
  \end{align} for $j\ge i\ge 1.$

  Write simply  $\xi_n:=\sum_{j=1}^n\eta_j$ for $n\ge1.$ Since $\xi_n$ is nondecreasing in $n,$ the limit $\xi:=\lim_{n\rto}\xi_n$ exists almost surely. Moreover, by L\'evy's monotone convergence theorem, we have
  \begin{align}\label{exm1}
    E\xi&=E\Bigg(\lim_{n\rto}\sum_{j=1}^n\eta_j\Bigg)=\sum_{j=1}^\infty P(\eta_j=1)=\sum_{j=1}^\infty\frac{1}{D(j)}=\zeta(D)<\infty,
  \end{align}
 where $\zeta(D)$ is defined in (\ref{zed}). Consequently, $\xi<\infty$ almost surely.

In view of \eqref{exs}, \eqref{depsi}, and \eqref{dns}, the $k$-th moment of $\xi_n$ is given by
  \begin{align}
  E\Bigg(\sum_{j=1}^n\eta_j\Bigg)^k=\sum_{m=1}^k\sum_{\begin{subarray}{c}l_1+\dots+l_m=k,\\
  l_s\ge1,\,s=1,...,m
  \end{subarray}}\frac{k!}{l_1!\cdots l_m!}\Psi_n(m),\no
  \end{align}
    where
  \begin{align*}
 \Psi_n(m)= \sum_{1\le j_1<...<j_m\le n}\frac{1}{D(j_1)D(j_2-j_1)\cdots D(j_{m}-j_{m-1})}.
\end{align*}
By Part (i) of Theorem \ref{prpd} and \eqref{exm1}, we  obtain
 \begin{align}\label{eexk}
  E\z(\xi^k\y)&=\lim_{n\rto} E\Bigg(\sum_{j=1}^n\eta_j\Bigg)^k\\
  &=\lim_{n\rto}\sum_{m=1}^k\sum_{\begin{subarray}{c}l_1+\dots+l_m=k,\\
  l_s\ge1,\,s=1,...,m
  \end{subarray}}\frac{k!}{l_1!\cdots l_m!}\sum_{1\le j_1<...<j_m\le n}\Psi_n(m)\no\\
  &=\sum_{m=1}^k\sum_{\begin{subarray}{c}l_1+\dots+l_m=k,\\
  l_s\ge1,\,s=1,...,m
  \end{subarray}}\frac{k!}{l_1!\cdots l_m!}\zeta(D)^m=\sum_{m=1}^k\sum_{\begin{subarray}{c}l_1+\dots+l_m=k,\\
  l_s\ge1,\, s=1,...,m
  \end{subarray}}\frac{k!}{l_1!\cdots l_m!} (E\xi)^m.\no
  \end{align}
For $k\ge2,$ by expanding recursively, we get
\begin{align}
  E(\xi^k)&=E\xi+\sum_{m=2}^k\sum_{\begin{subarray}{c}
l_1+...+l_m=k\\
l_{s}\ge 1,\,s=1,...,m
  \end{subarray}}\frac{k!}{l_1!l_{2}!\cdots l_{m}!}\z(E\xi\y)^m\no\\
  &=E\xi+\sum_{m=2}^k\sum_{l_1=1}^{k-m+1}\frac{k!}{l_1!\z(k-l_1\y)!}\sum_{\begin{subarray}{c}
l_2+...+l_m=k-l_1\\
l_{s}\ge 1,\,s=2,...,m
  \end{subarray}}\frac{\z(k-l_1\y)!}{l_{2}!\cdots l_{m}!}\z(E\xi\y)^m\no\\
  &=E\xi+\sum_{l_1=1}^{k-1}\frac{k!}{l_1!\z(k-l_1\y)!}\sum_{m=2}^{k-l_1+1}\sum_{\begin{subarray}{c}
l_2+...+l_m=k-l_1\\
l_{s}\ge 1,\,s=2,...,m
  \end{subarray}}\frac{\z(k-l_1\y)!}{l_{2}!\cdots l_{m}!}\z(E\xi\y)^m\no\\
&=E\xi+E\xi\sum_{l_1=1}^{k-1}\frac{k!}{l_1!\z(k-l_1\y)!}\sum_{m=1}^{k-l_1}\sum_{\begin{subarray}{c}
h_1+...+h_m=k-l_1\\
h_{s}\ge 1,\,s=1,...,m
  \end{subarray}}\frac{\z(k-l_1\y)!}{h_{1}!\cdots h_{m}!}\z(E\xi\y)^m\cr
&=E\xi+E\xi\sum_{l_1=1}^{k-1}\frac{k!}{l_1!\z(k-l_1\y)!}E\z(\xi^{k-l_1}\y).\no
\end{align}
Here the last equality follows from (\ref{eexk}). Noticing that $E(\xi^0)=1,$  then consulting to \eqref{exm1} and \eqref{eexk}, we conclude that
$\xi $ satisfies $E\xi=\zeta(D)$ and
\begin{align}
   E\z(\xi^k\y)&=E\xi+E\xi\sum_{j=1}^{k-1}\frac{k!}{j!(k-j)!}E\z(\xi^{k-j}\y)
   =E\xi \z(E\z((\xi+1)^k\y)- E\z(\xi^k\y)\y).\no
\end{align}
Consequently, we get
\begin{align}
    \frac{E\z((\xi+1)^k\y)}{E(\xi^k)}= \frac{1+E(\xi)}{E(\xi)}\text{ for } k\ge1 \text{ and } E\xi=  \zeta(D).\no
\end{align}
  Suppose $\eta$ is a random variable and $\eta\sim \mathrm{Geo}\z(\frac{1}{\zeta(D)+1}\y).$
Then  direct computation shows that
\begin{align}
  \frac{E\z((\eta+1)^k\y)}{E(\eta^k)}=\frac{1+E\eta}{E\eta}\text{ for } k\ge 1\text{ and } E\eta=\zeta(D).\no
\end{align}
Thus, $\xi$ and $\eta$ share identical moments. Since there exists $C>0$ such that for $|t|<C$, $$E\z(e^{t\eta}\y)=\frac{\z(1+\zeta(D)\y)^{-1}}{1-e^t\zeta(D)\z(1-\zeta(D)\y)^{-1}}<\infty,$$
according to \cite[Theorem 30.1]{Billingsley}, the moment sequence uniquely determines the distribution.  Therefore,
\[\xi\sim \mathrm{Geo}\z(\frac{1}{\zeta(D)+1}\y).\]
This completes the proof of Part (i) of Theorem \ref{thbb}. \qed

  \subsection{Proof of Part (ii) of Theorem \ref{thbb}}

  \proof
   Assume $0\le \sigma\le 1$. Fix $0<\ve<1.$ Let $n_1$ be as in Part (ii) of Theorem \ref{thbb}. By assumption, we have
  \begin{align}\label{dnul}
    1-\ve\le \frac{\lambda_\sigma}{\rho(0,n)}\le 1+\ve,\quad    1-\ve\le \frac{\lambda_\sigma D(n)}{\rho(i,i+n)}\le 1+\ve, \quad n\ge n_1,\ i\ge n_1.
  \end{align} In order to apply Theorem \ref{prpd}, we set $n_0=n_1.$
  Let $j_0=0$. For $m\geq1$, put
 \begin{align}
    &H=\{(j_1,...,j_m)\mid 1\le j_1<....<j_m\le n,\ j_{i+1}-j_i>n_0,\  \forall~0\le i\le m-1\},\label{dsa}\\
    &A=\{(j_1,...,j_m)\mid 1\le j_1<....<j_m\le n,\ (j_1,...,j_m)\notin H\}.\label{dss}
       \end{align}
Fix an integer $k\ge 1.$ From \eqref{exs}, for $1\le m\le k$,
  \begin{align}\label{eks}
  E&\Bigg(\sum_{j=1}^n\eta_j\Bigg)^k=\sum_{m=1}^k\sum_{\begin{subarray}{c}l_1+\dots+l_m=k,\\
  l_s\ge1,\,s=1,...,m
  \end{subarray}}\frac{k!}{l_1!\cdots l_m!}\Psi_n(m),
  \end{align} where
 $$
    \Psi_n(m)=\sum_{1\le j_1<...<j_m\le n}\frac{1}{\rho(0,j_1)\rho(j_1,j_2)\cdots \rho(j_{m-1},j_m)}.
 $$
Next, we aim to prove
\begin{align}\label{lpa}
  \lim_{n\rto}\frac{\Psi_n(m)}{\z(S(n)\y)^m}=\lambda_\sigma^{-m}.
\end{align}
 For a set $V\subset\mathbb{N}^m,$ write
 \begin{align}\label{dpc}
   \Psi_n(m,V)=\sum_{(j_1,...,j_m)\in V} \frac{1}{\rho(0,j_1)\rho(j_1,j_2)\cdots \rho(j_{m-1},j_m)}.
 \end{align}
Then,
 \begin{align}\label{gab}
  \Psi_n(m)=\Psi_n(m,H)+\Psi_n(m,A),
  \end{align} where $H$ and $A$ are defined in \eqref{dsa} and \eqref{dss}.

  Consider first the term $\Psi_n(m,H).$ From \eqref{dnul}, we get
  \begin{align*}
(1-\ve)^m&\lambda_\sigma^{-1}\sum_{(j_1,...,j_m)\in H} \frac{1}{D(j_1)D(j_2-j_1)\cdots D(j_m-j_{m-1})}\cr
   &\le\Psi_n(m,H)=\sum_{(j_1,...,j_m)\in H}\frac{1}{\rho(0,j_1)\rho(j_1,j_2)\cdots \rho(j_{m-1},j_m)}\\
     &\le(1+\ve)^m\lambda_\sigma^{-1} \sum_{(j_1,...,j_m)\in H} \frac{1}{D(j_1)D(j_2-j_1)\cdots D(j_m-j_{m-1})}.\no
\end{align*}
Then applying Theorem  \ref{prpd},  we obtain
  \begin{align*}
  (1-\ve)^m\lambda_\sigma^{-m}\le \varliminf_{n\rto} \frac{\Psi_n(m,H)}{\z(S(n)\y)^{m}} \le \varlimsup_{n\rto} \frac{\Psi_n(m,H)}{\z(S(n)\y)^{m}}\le (1+\ve)^m\lambda_\sigma^{-m}.
  \end{align*}
  Taking $\ve\to 0,$ we get
  \begin{align}\label{lpn}
  \lim_{n\rto}\frac{\Psi_n(m,H)}{\z(S(n)\y)^{m}}=\lambda_\sigma^{-m}.
\end{align}
  Next, we consider  the term $\Psi_n(m,A)$ in \eqref{gab}. We shall show that
  \begin{align}
    \lim_{n\rto}\frac{\Psi_n(m,A)}{\z(S(n)\y)^{m}}=0.\label{gbi}
  \end{align}
  To this end, for
  $0\le i_1<...<i_l\le m-1$ and $k_1,...,k_l\ge1,$  set
  $$A\z(i_1,k_1,...,i_l,k_l\y)=\left\{(j_1,...,j_m)\left|\begin{array}{c}1\le j_1<....<j_m\le n, \\
    j_{i_s+1}-j_{i_s}=k_s,\, 0\le s\le l,\\
    j_{i+1}-j_i>n_0,\, i\in\{0,..,m-1\}/\{i_1,...,i_l\}\end{array}\right.\right\}.$$
  Then, we have
  \begin{align*}
    A=\bigcup_{l=1}^{m-1}\bigcup_{\begin{subarray}{c}
                                    0\le i_1<...<i_l\le m-1\\
                                     1\le k_1,...,k_l\le n_0
                                  \end{subarray}
    }A\z(i_1,k_1,...,i_l,k_l\y).
  \end{align*}
  Thus, it follows that
  \begin{align}
    \Psi_n(m,A)=\sum_{l=1}^{m-1}\sum_{\begin{subarray}{c}
                                    0\le i_1<...<i_l\le m-1\\
                                     1\le k_1,...,k_l\le n_0
                                  \end{subarray}}
    \Psi_n\z(m,{A(i_1,k_1,...,i_l,k_l)}\y).\label{gbd}
  \end{align}
 Fix $1\le l\le m,$ $0\le i_1<...<i_l\le m-1$ and $1\le k_1,...,k_l\le n_0.$ By \eqref{gbd}, to prove \eqref{gbi}, it suffices to show that
  \begin{align}
   \lim_{n\rto}\frac{1}{(S(n))^{m}} \Psi_n\z(m,{A(i_1,k_1,...,i_l,k_l)}\y)=0.\label{bik0}
  \end{align}
In fact, by \eqref{dnul} and the fact $\rho(i,j)\ge 1,$ we have
    \begin{align}
    & \Psi_n\z(m,{A(i_1,k_1,...,i_l,k_l)}\y)=\sum_{(j_1,...,j_m)\in A\z(i_1,k_1,...,i_l,k_l\y)} \frac{1}{\rho\z(0,j_1\y)\rho\z(j_1,j_2\y)\cdots \rho\z(j_{m-1},j_m\y)}\no\\
    &\le\sum_{(j_1,...,j_m)\in A(i_1,k_1,...,i_l,k_l)} \;\prod_{i\in \{0,..,m-1\}/\{i_1,...,i_l\}}\frac{1}{\rho\z(j_i, j_{i+1}\y)}\no\\
    &\le (\lambda_{\sigma}^{-1}\vee 1)\sum_{(j_1,...,j_m)\in A(i_1,k_1,...,i_l,k_l)} \;\prod_{i\in \{0,..,m-1\}/\{i_1,...,i_l\}}\frac{1+\ve}{D\z(j_{i+1}-j_i\y)}\no\\
    & =(\lambda_{\sigma}^{-1}\vee 1)\sum_{1\le j_1'<...<j'_{m-l}\le n-\sum_{s=1}^lk_s}(1+\varepsilon)^{m-l}\frac{1}{D\z(j_1'\y)D\z(j_2'-j_1'\y)\cdots D\z(j_{m-l}'-j_{m-l-1}'\y)}.\no
  \end{align}
   Therefore, we conclude that
\begin{align}
\varlimsup_{n\rto}\frac{1}{(S(n))^{m}}&\Psi_n\z(m,{A(i_1,k_1,...,i_l,k_l)}\y)
   \le \varlimsup_{n\rto}\frac{ (1+\varepsilon)^{m-l}}{(S(n))^{l}}\frac{\lambda_{\sigma}^{-1}\vee 1 }{(S(n))^{m-l}}\no\\
   &\times\sum_{1\le j_1'<...<j'_{m-l}\le n-\sum_{s=1}^lk_s}\frac{1}{D\z(j_1'\y)D\z(j_2'-j_1'\y)\cdots D\z(j_{m-l}'-j_{m-l-1}'\y)}=0,\no
  \end{align}
where we use Theorem \ref{prpd} to get the last step.
Thus, \eqref{bik0} is true and so is \eqref{gbi}.
\par
  Dividing by $(S(n))^{m}$ on both sides of \eqref{gab} and taking $n\to\infty,$ owing to \eqref{lpn} and \eqref{gbi}, we get \eqref{lpa}. Now, putting \eqref{eks} and \eqref{lpa}  together, we deduce that
\begin{align}
  \lim_{n\rto}\frac{E\z(\sum_{j=1}^n\eta_j\y)^k}{(S(n))^{k}}=k!\lambda_\sigma^{-k}.\no
\end{align}
Finally, since $ [(2n)!]^{-\frac{1}{2n}}\sim \frac{e }{2n}$ as $n\rto,$ we get
$\sum_{k=1}^\infty \z[(2k)!\lambda_\sigma^{-2k}\y]^{-\frac{1}{2k}}=\infty.$
By Carleman's criterion (see e. g., \cite[Chap. II, \S 12]{s}), the moment problem has a unique solution. Hence,
\begin{align*}
  \frac{\sum_{j=1}^n\eta_j}{S(n)}\overset{d}\to \xi,\text{ as }n\rto,
\end{align*}
   where $\xi\sim \mathrm{Exp}\z(\lambda_\sigma\y)$. Part (ii) of Theorem \ref{thbb} is proved. \qed
\subsection{Proof of Theorem \ref{tha}}
\proof Based on Proposition \ref{thg}, the proof of Theorem \ref{tha} parallels Part (ii) of Theorem \ref{thbb}. We outline the argument below for completeness.

   Assume $\alpha>0.$ Fix $0<\ve<1.$ Let $n_2$ be the one in Theorem \ref{tha}. From \eqref{al}, we have
  \begin{align}\label{dnul2}
    1-\ve\le \frac{\beta n}{\rho(0,n)}\le 1+\ve, \quad  1-\ve\le \frac{\beta (i+n)^{1-\alpha}((i+n)^\alpha-i^\alpha)}{\rho(i,i+n)}\le 1+\ve, \quad n\ge n_2,\ i\ge n_2.
  \end{align}
   Fix $k\ge 1.$ From \eqref{exs}, we see that
  \begin{align}\label{eks2}
  E&\Bigg(\sum_{j=1}^n\eta_j\Bigg)^k=\sum_{m=1}^k\sum_{\begin{subarray}{c}l_1+\dots+l_m=k,\\
  l_s\ge1,\, s=1,...,m
  \end{subarray}}\frac{k!}{l_1!\cdots l_m!}\Psi_n(m),
  \end{align} where
 $$
    \Psi_n(m)=\sum_{1\le j_1<...<j_m\le n}\frac{1}{\rho(0,j_1)\rho(j_1,j_2)\cdots \rho(j_{m-1},j_m)}.
 $$
Next, we  prove
\begin{align}\label{lpa2}
  \lim_{n\rto}\frac{\Psi_n(m)}{(\log n)^m}=\frac{\prod_{j=0}^{m-1}(j+\alpha)}{m!\alpha^m \beta^m}, \quad m\ge1.
\end{align}
To this end, let $j_0=0$ and set $n_0=n_2.$
 For $m\ge 1,$ let $H$ and $A$ be the sets defined in \eqref{dsa} and \eqref{dss}. Then  we have
 \begin{align}\label{gab2}
  \Psi_n(m)=\Psi_n(m,H)+\Psi_n(m,A),
  \end{align} where for a set $V\subset\mathbb{N}^m,$ $\Psi_n(m,V)$ is defined by \eqref{dpc}.
  Due to \eqref{dnul2}, we have
  \begin{align*}
(1-\ve)^m\beta^{-m}&\sum_{(j_1,...,j_m)\in H} \frac{1}{j_1j_2^{1-\alpha}(j_2^{\alpha}-j_1^{\alpha})\cdots j_m^{1-\alpha}(j_m^{\alpha}-j_{m-1}^{\alpha})}\\
   &\le\Psi_n(m,H)=\sum_{(j_1,...,j_m)\in H}\frac{1}{\rho(0,j_1)\rho(j_1,j_2)\cdots \rho(j_{m-1},j_m)}\\
     &\le(1+\ve)^m\beta^{-m}\sum_{(j_1,...,j_m)\in H} \frac{1}{j_1j_2^{1-\alpha}(j_2^{\alpha}-j_1^{\alpha})\cdots j_m^{1-\alpha}(j_m^{\alpha}-j_{m-1}^{\alpha})}.\no
\end{align*}
Applying Proposition \ref{thg} with $b_s\equiv n_0,\ s=1,...,m,$ we get
  \begin{align}\label{gau2}
  (1-\ve)^m\frac{\prod_{j=0}^{m-1}(j+\alpha)}{m!\alpha^m\beta^{m}}&\le \varliminf_{n\rto} \frac{\Psi_n(m,H)}{(\log n)^{m}} \\
  &\le \varlimsup_{n\rto} \frac{\Psi_n(m,H)}{(\log n)^{m}}\le (1+\ve)^m \frac{\prod_{j=0}^{m-1}(j+\alpha)}{m!\alpha^m\beta^{m}}.\no
  \end{align}Since $\ve$ is arbitrary, letting $\ve\to 0$ in \eqref{gau2}, we deduce that
  \begin{align}\label{lpn2}
  \lim_{n\rto}\frac{\Psi_n(m,H)}{(\log n)^{m}}=\frac{\prod_{j=0}^{m-1}(j+\alpha)}{m!\alpha^m\beta^{m}}.
\end{align}
Using \eqref{ijlu} and Proposition \ref{thg}, analogous to \eqref{gbi}, we  show that
  \begin{align}
    \lim_{n\rto}\frac{\Psi_n(m,A)}{(\log n)^{m}}=0.\label{gbi2}
  \end{align}
      Putting \eqref{gab2}, \eqref{lpn2} and \eqref{gbi2} together, we get \eqref{lpa2}. With \eqref{lpa2} in hand, dividing \eqref{eks2} by $(\log n)^k$ and taking $n\to\infty$, we infer that
\begin{align}
  \lim_{n\rto}\frac{E\z(\sum_{j=1}^n\eta_j\y)^k}{(\log n)^{k}}=\frac{\prod_{j=0}^{k-1}(j+\alpha)}{\alpha^k\beta^{k}}.\no
\end{align}
Let  $\xi\sim \mathrm{Gamma}(\alpha,1)$. Then
  \begin{align}
    E\z(\xi^k\y)&=\frac{1}{\Gamma(\alpha)}\int_{0}^\infty x^{k+\alpha-1}e^{-x}dx=\frac{\Gamma(k+\alpha)}{\Gamma(\alpha)}=\prod_{j=0}^{k-1}(j+\alpha).\no
  \end{align}
 Using Stirling's approximation, we get
 \[
\z[E\z({\xi}^{2n}\y)\y]^{-\frac{1}{2n}}=\z(\prod_{j=0}^{2n-1}(j+\alpha)\y)^{-\frac{1}{2n}}\ge (\lfloor 2n+\alpha\rfloor !)^{-\frac{1}{2n}}\sim \frac{e}{2n}\]
 as $n\rto.$ Thus
 $\sum_{n=1}^\infty\z(E\z(\xi^{2n}\y)\y)^{-\frac{1}{2n}}=\infty.$
Applying Carleman's test for the uniqueness of the moment problem (see e. g., \cite[Chap.II, \S 12]{s}), we have
\begin{align*}
  \frac{\alpha \beta\sum_{j=1}^n\eta_j}{\log n}\overset{d}\to \xi, \text{ as } n\rto,
\end{align*}
  where $\xi\sim \mathrm{Gamma}(\alpha, 1)$.  This completes the proof of Theorem \ref{tha}. \qed

\section{Applications in branching processes}\label{sbp}
In this section,  we give the proofs of Corollaries \ref{c3}--\ref{thz}.

\subsection{Proof of Corollary \ref{c3}}
\proof As usual we use $\mathbb P_x$ and $\mathbb E_x$ to denote the probability measure and the corresponding expectation operator  induced by a $d$-dimensional Brownian motion starting at $x\in\mathbb R^d.$
For $r>0,$
define the stoping time
$$T_r=\inf\{t>0:B_t\in \mathcal S(0,r)\}.$$
Fix $0<r<R<\infty$ and $x\in\mathbb R^d$ such that $r<|x|<R.$ Then by solving the Dirichlet problem on the annulus $A=\{y\in \mathbb R^d: r<|y|<R\},$ we obtain(see \cite[Theorem 3.18]{mp})
\begin{align}\label{prr}
  \mathbb P_x(T_r>T_R)=\frac{|x|^{2-d}-r^{2-d}}{R^{2-d}-r^{2-d}}.
\end{align}
Define $\eta_j=1_{\{j\in C(a,b)\}}$ for $j\in \mathbb N,$ so that  $$|C(a,b)\cap [1,n]|=\sum_{j=1}^n\eta_j,\quad n\ge1.$$
It follows immediately from \eqref{prr} that
\begin{align*}
  \mathbb P_0(\eta_i=1)=\frac{(bi)^{2-d}-(bi+a)^{2-d}}{(bi)^{2-d}}, \ i\ge 1.
\end{align*}
Now fix $m\ge 1$ and $1\le j_1<j_2<...<j_m<\infty.$ Using  \eqref{prr} and the strong Markov property, we infer that
\begin{align}
  &\mathbb P_0\z(\eta_{j_1}=1,...,\eta_{j_m}=1\y)=\mathbb P_0\z(j_1\in C(a,b),...,j_m\in C(a,b)\y)\no\\
  &=\frac{(j_1b+a)^{2-d}-(j_1b)^{2-d}}{(j_2b+a)^{2-d}-(j_1b)^{2-d}}\cdots \frac{(j_{m-1}b+a)^{2-d}-(j_{m-1}b)^{2-d}}{(j_mb+a)^{2-d}-(j_{m-1}b)^{2-d}} \frac{(j_{m}b)^{2-d}-(j_{m}b+a)^{2-d}}{(j_mb)^{2-d}}\no\\
  &=\frac{(j_1+a/b)^{2-d}-j_1^{2-d}}{(j_2+a/b)^{2-d}-j_1^{2-d}}\cdots \frac{(j_{m-1}+a/b)^{2-d}-j_{m-1}^{2-d}}{(j_m+a/b)^{2-d}-j_{m-1}^{2-d}} \frac{j_{m}^{2-d}-(j_{m}+a/b)^{2-d}}{j_m^{2-d}}\no.
\end{align}
Similar to \eqref{exs}, we have
\begin{align}\label{eek}
  \mathbb E_0\Bigg(\sum_{j=1}^n\eta_j\Bigg)^k&=\sum_{m=1}^k\sum_{\begin{subarray}{c}l_1+\dots+l_m=k,\\
  l_s\ge1,\,s=1,...,m
  \end{subarray}}\frac{k!}{l_1!\cdots l_m!}\Psi_n(m),\quad k\ge 1.
  \end{align}
  where for $n\ge m \ge 1,$
  \begin{align*}
    \Psi_n(m)=\sum_{1\le j_1<...<j_m\le n}&\frac{(j_1+a/b)^{2-d}-j_1^{2-d}}{(j_2+a/b)^{2-d}-j_1^{2-d}}\cdots \frac{(j_{m-1}+a/b)^{2-d}-j_{m-1}^{2-d}}{(j_m+a/b)^{2-d}-j_{m-1}^{2-d}} \\
    &\times\frac{j_{m}^{2-d}-(j_{m}+a/b)^{2-d}}{j_m^{2-d}}.
  \end{align*}
  Clearly, we have
  \begin{align}
    &(j+a/b)^{2-d}-j^{2-d}\sim a/b(2-d)j^{1-d},\quad j\to\infty.\label{sjab}\end{align}
   Moreover,  by the mean value theorem, there exist $\kappa_1\in (j,j+a/b),$  $\kappa_2\in (i,j)$ such that
  \begin{align*}
    \frac{(j+a/b)^{2-d}-i^{2-d}}{j^{2-d}-i^{2-d}}=1+\frac{(j+a/b)^{2-d}-j^{2-d}}{j^{2-d}-i^{2-d}}
    =1+\frac{a/b(2-d)\kappa_1^{1-d}}{(2-d)(j-i)\kappa_2^{1-d}},
  \end{align*}
  which implies
  \begin{align}
    1<\frac{(j+a/b)^{2-d}-i^{2-d}}{j^{2-d}-i^{2-d}}<1+\frac{a}{b(j-i)}, \quad j>i\ge 1.\label{jiab}
  \end{align}
  Fix $\ve>0.$ Due to \eqref{sjab} and \eqref{jiab}, there exists $n_3>0$ such that
  \begin{gather}\label{jul}
   a/b(d-2)j^{1-d}(1-\ve)\le j^{2-d}-(j+a/b)^{2-d}\le a/b(d-2)j^{1-d}(1+\ve),\quad j\ge n_3,\\
   i^{2-d}-j^{2-d}<i^{2-d}-(j+a/b)^{2-d}<(1+\ve)(i^{2-d}-j^{2-d}),\quad j-i\ge n_3.\label{jiul}
  \end{gather}
    Let $n_0=n_3$ and define $H$ and $A$ as in \eqref{dsa} and \eqref{dss}.
    Then,
 \begin{align}\label{gaba}
  \Psi_n(m)=\Psi_n(m,H)+\Psi_n(m,A),
  \end{align} where
     \begin{align*}
   \Psi_n(m,V)=\sum_{(j_1,...,j_m)\in V}& \frac{(j_1+a/b)^{2-d}-j_1^{2-d}}{(j_2+a/b)^{2-d}-j_1^{2-d}}\cdots \frac{(j_{m-1}+a/b)^{2-d}-j_{m-1}^{2-d}}{(j_m+a/b)^{2-d}-j_{m-1}^{2-d}} \\&\times \frac{j_{m}^{2-d}-(j_{m}+a/b)^{2-d}}{j_m^{2-d}}
 \end{align*}
for $V=A$ or $H.$
Owing to \eqref{jul} and \eqref{jiul}, we have
\begin{align}\label{pmhu}
  \Psi_n(m,H)&\le (a/b(d-2)(1+\ve))^{m}\sum_{(j_1,...,j_m)\in H} \z(\prod_{l=1}^{m-1}\frac{j_l^{1-d}}{j_{l}^{2-d}-j_{l+1}^{2-d}}\y)\frac{1}{j_m}\\
  &=(a/b(d-2)(1+\ve))^{m}\sum_{(j_1,...,j_m)\in H} \frac{1}{j_1}\prod_{l=1}^{m-1}\frac{1}{j_l^{3-d}\z(j_{l+1}^{d-2}-j_{l}^{d-2}\y)}\no
\end{align}
and
\begin{align}\label{pmhl}
\Psi_n(m,H)&\ge \z(\frac{a/b(d-2)(1-\ve)}{1+\ve}\y)^{m}\sum_{(j_1,...,j_m)\in H} \z(\prod_{l=1}^{m-1}\frac{j_l^{1-d}}{j_{l}^{2-d}-j_{l+1}^{2-d}}\y)\frac{1}{j_m}\\
  &=\z(\frac{a/b(d-2)(1-\ve)}{1+\ve}\y)^{m}\sum_{(j_1,...,j_m)\in H} \frac{1}{j_1}\prod_{l=1}^{m-1}\frac{1}{j_l^{3-d}\z(j_{l+1}^{d-2}-j_{l}^{d-2}\y)}.\no
\end{align}
Based on \eqref{pmhu} and \eqref{pmhl}, putting $\alpha=d-2,$ $b_s=n_0,$  $s=1,2,...,m$ and applying Proposition \ref{thg}, we get
\begin{align*}
  \z(\frac{a/b(d-2)(1-\ve)}{1+\ve}\y)^{m}&\frac{\prod_{j=0}^{m-1}(j+d-2)}{m!(d-2)^m}\le \varliminf_{n\rto} \frac{\Psi_n(m,H)}{\z(\log n\y)^{m}}\no\\
   &\le \varlimsup_{n\rto} \frac{\Psi_n(m,H)}{\z(\log n\y)^{m}}\le (a/b(d-2)(1+\ve))^{m}\frac{\prod_{j=0}^{m-1}(j+d-2)}{m!(d-2)^m}.
  \end{align*}
Since $\ve$ is arbitrary, we have
\begin{align}\label{ph}
  \lim_{n\rto} \frac{\Psi_n(m,H)}{\z(\log n\y)^{m}}= (a/b)^{m}\frac{\prod_{j=0}^{m-1}(j+d-2)}{m!}.
  \end{align}
  On the other hand, based on \eqref{jul} and \eqref{jiul}, by an argument similar to the proof of \eqref{gbd}, we show that
   \begin{align}\label{pa}
     \lim_{n\rto} \frac{\Psi_n(m,A)}{\z(\log n\y)^{m}}=0.
   \end{align}
  In view of \eqref{gaba}, \eqref{ph} and \eqref{pa}, dividing \eqref{eek} by $(a/b\log n)^k$ and taking the limit, we obtain
  \begin{align*}
    \lim_{n\rto}\frac{1}{(a/b\log n)^k}\mathbb E_0\Bigg(\sum_{j=1}^n\eta_j\Bigg)^k=\prod_{j=0}^{k-1}(j+d-2), \quad k\ge1.
  \end{align*}
  Applying again Carleman's test for the uniqueness of the moment problem (see e. g., \cite[Chap. II, \S 12]{s}), we have
\begin{align*}
  \frac{ \sum_{j=1}^n\eta_j}{a/b\log n}\overset{d}\to \xi, \text{ as } n\rto,
\end{align*}
  where $\xi\sim \mathrm{Gamma}(d-2, 1)$.  This completes the proof of Corollary  \ref{c3}. \qed

\subsection{Proof of Corollary \ref{c4}}
\proof To prove Corollary \ref{c4}, recall that the generator of the geometric Brownian motion $\{X_t\}_{t\ge0}$ defined in \eqref{sde} or \eqref{geob} is given by
$$Lf(x)=\frac{1}{2}\sigma^2x^2f''(x)+\mu xf'(x),$$ where $f$ is bounded and twice continuously differentiable. A non-degenerate  solution $S(x)$ of the equation    \begin{align}
  \frac{1}{2}\sigma^2x^2S''(x)+\mu xS'(x)=0, \text{ or }LS=0\label{hf}
\end{align}
is called the scale function of $\{X_t\}_{t\ge 0}.$
 It is easy to check that \begin{align}\label{sx}
  S(x)=\z\{\begin{array}{ll}
 x^{1-2\mu/\sigma^2}, & \text{if } 2\mu/\sigma^2\ne 1,\\
 \log x, &\text{if } 2\mu/\sigma^2=1
  \end{array}\y.
\end{align}
is a special solution of equation \eqref{hf}.

In what follows, we denote by $P_x$ and $E_x$ the probability measure and the corresponding expectation operator of a geometric Brownian motion starting from $x.$
For $y>0,$ define $$T_y=\inf\{t>0: X_t=y \}.$$
Then, for $0<r<x<R<\infty$ we have (see for example \cite[Theorem 6.12]{fck})
\begin{align}\label{exgb}
  P_x(T_r>T_R)=\frac{S(x)-S(r)}{S(R)-S(r)}.
\end{align}
If $2\mu/\sigma^2>1,$ then from \eqref{sx} and \eqref{exgb} we obtain
\begin{align}\label{exgbg}
  P_x(T_r>T_R)=\frac{x^{1-2\mu/\sigma^2}-r^{1-2\mu/\sigma^2}}{R^{1-2\mu/\sigma^2}-r^{1-2\mu/\sigma^2}},
\end{align}
for $ 0<r<x<R<\infty.$

Based on \eqref{exgbg}, the rest of the proof proceeds verbatim as the analogous part of the proof of Corollary \ref{c3} (after \eqref{prr}).  \qed

\subsection{Proof of Corollary \ref{thy}}
\proof Let $f_0(s)=s$ for $s\in[0,1]$. Recall that $f(s)=1/(2-s)$ for $s\in[0,1]$. Define recursively $f_n(s)=f(f_{n-1}(s))$ for $n\ge 1.$ By induction, we get   $f_n(s)=1-\frac{1-s}{1+n(1-s)}$ for $n\ge 0$  with derivatives  $f^{(j)}_n(0)=\frac{j!n^{j-1}}{(1+n)^{j+1}}$ for $j\ge 1.$ For the Galton-Watson process $\{Y_n\}_{n\geq0}$, the generating function satisfies (see e. g., \cite[p.2]{an})
\begin{align}
  E\z(s^{Y_{k+n}}\,\middle|\,Y_k=1\y)=f_n(s), \quad k\ge 0,\ n\ge 0. \no
\end{align}
Set $\eta_i=1_{\left\{Y_i=1\right\}}$ for  $i\ge 1.$ Then  $P(\eta_i=1)=P(Y_i=1)=f_i'(0)=\frac{1}{(1+i)^2},$ $i\ge1.$ By the Markov property of $\{Y_n\}_{n\ge 1},$ for $j\ge i\ge 1$, we have
\begin{align*}
  P\z(\eta_{j}=1\,\middle|\,\eta_i=1,...,\eta_1=1\y)&=P\z(Y_{j}=1\,\middle|\,Y_i=1\y)=f_{j-i}'(0)=\frac{1}{(j-i+1)^2}.
\end{align*}
Therefore, \eqref{mc} holds with $D(n)=(1+n)^2,$ $n\ge1$ and $\rho(i,j)=D(j-i)=(j-i+1)^2,\ j\ge i\ge 0.$
Notice that
\begin{align*}
  \zeta(D):=\sum_{n=1}^\infty \frac{1}{D(n)}=\sum_{n=1}^\infty\frac{1}{(1+n)^2}=\frac{\pi^2}{6}-1.
\end{align*}
Thus, an application of Part (i) of Theorem \ref{thbb} yields that
\[\sum_{i=1}^n\eta_i\overset{a.s.}{\to} \xi\text{ as }n\rto,\] where $\xi\sim \mathrm{Geo}(6/\pi^2).$
 Corollary \ref{thy} is proved. \qed

\subsection{Proof of Corollary \ref{thz}}
\proof Recall that $\{Z_n\}_{n\ge 0}$ is a BPVE with immigration which satisfies $P(Z_0=0)=1$ and
\begin{align}\label{mb2}
  E\z(s^{Z_{k+1}}\,\middle|\,Z_{k},..,Z_0\y)=\z(\tilde f_{k+1}(s)\y)^{Z_{k}+1}, \quad k\ge 0,
\end{align} where $\tilde f_k,\ k\ge 1$ are given in \eqref{fis}.
  For $j\ge i\ge 0$, define
 \begin{align*}
   F_{i,j}(s)=E\z(s^{Z_{j}}\,\middle |\, Z_i=0\y),\quad s\in [0,1].
 \end{align*}
 By induction and \eqref{mb2}, for $j\ge i\ge 0$, we have
\begin{align}\label{fn}
  F_{i,j}(s)=\prod_{k=i+1}^j \tilde f_{k,j}(s),\quad s\in [0,1],
\end{align}
where
$\tilde f_{k,j}(s)=\tilde f_k\z(\tilde f_{k+1}\z(\cdots\z(\tilde f_j(s)\y)\y)\y)$ for  $j\ge k\ge 1.$
Recall that $m_k=\tilde f_k'(1)=\frac{1-p_k}{p_k}$ for $k\ge1.$ Thus some direct calculation yields
$\tilde f_k(s)=1-\frac{m_k(1-s)}{1+m_k(1-s)},\ k\ge1.$ As a consequence, by induction we deduce that
\begin{align}
  \tilde f_{k,j}(s)&=1-\frac{m_k\cdots m_j(1-s)}{1+\sum_{t=k}^j m_t\cdots m_j(1-s)}
  =\frac{1+\sum_{t=k+1}^j m_t\cdots m_j(1-s)}{1+\sum_{t=k}^j m_t\cdots m_j(1-s)},\quad j\ge k\ge1.\label{fkn}
\end{align}
Substituting \eqref{fkn} into \eqref{fn}, we obtain
\begin{align*}
  F_{i,j}(s)=\frac{1}{1+\sum_{t=i+1}^j m_t\cdots m_j(1-s)},\quad j\ge i\ge 0.
\end{align*}
For $j\ge i\ge 0,$ writing \begin{align*}
  \rho(i,j)=1+\sum_{t=i+1}^j m_t\cdots m_j,
\end{align*} then
\begin{align}\no
  P\z(Z_j=0\,\middle|\,Z_i=0\y)=F_{i,j}(0)=\frac{1}{1+\sum_{t=i+1}^j m_t\cdots m_j}=\frac{1}{\rho(i,j)}.
\end{align}
Let $\eta_i=1_{\{Z_i=0\}}$ for $i\ge 1.$ Then, for $j\ge i\geq1$, we have
\begin{align}
&P(\eta_i=1)=P(Z_i=0)=F_{0,i}(0)=\frac{1}{\rho(0,i)},\label{pz0}\\
 & P(\eta_j=1\,|\,\eta_i=1,...,\eta_0=1)=P(Z_j=0\, |\,Z_i=0)=\frac{1}{\rho(i,j)}. \label{zdij}
\end{align}

 \begin{lemma} \label{lmmd}
    Let $p_i,\ i\ge1$ be as in \eqref{dpa}. (i) Suppose $\sum_{k=i}^\infty r_k \to 0$ as $i \to \infty$. Then for any $\varepsilon > 0$, there exists $N \ge 1$ such that for all $i \ge N$ and $j - i > N$,
\begin{align}\label{d0n}
    1 - \varepsilon \le \frac{i}{\rho(0,i)} \le 1 + \varepsilon, \quad
    1 - \varepsilon \le \frac{j-i}{\rho(i,j)} \le 1 + \varepsilon.
\end{align}

(ii) Fix $B \in [0,1)$. If $r_i = \frac{B}{i}$ for $i \ge 1$, then for any $\varepsilon > 0$, there exists $N \ge 1$ such that for all $i \ge N$ and $j - i \ge N$,
\begin{align}\label{d1n}
    1 - \varepsilon \le \frac{i}{(1-B)\rho(0,i)} \le 1 + \varepsilon, \quad
    1 - \varepsilon \le \frac{j^{B}\z(j^{1-B} - i^{1-B}\y)}{(1-B)\rho(i,j)} \le 1 + \varepsilon.
\end{align}
 \end{lemma}
 \proof Observe that as $k\rto,$
\begin{align}\label{567hyhyt1}
m_k=\frac{1-p_k}{p_k}=\frac{1/2+r_k/4}{1/2-r_k/4}=1+r_k+O\z(r_k^{2}\y).
\end{align}
By Taylor expansion, there exists $\theta_k\in (0, m_k-1)$ such that
\begin{align*}
  \log m_k=m_k-1-\frac{(m_k-1)^2}{2(1+\theta_k)^2}.
\end{align*}
Combining this with \eqref{567hyhyt1}, we obtain for $j\ge t-1\ge 0,$
\begin{align}\label{46hyyr4}
m_t \cdots m_j&= \exp\left(\sum_{k=t}^j \log m_k \right)
= \exp\left(\sum_{k=t}^j \left[m_k - 1 - \frac{(m_k - 1)^2}{2(1 + \theta_k)^2}\right]\right)\\
&= \exp\left(\sum_{k=t}^j \left[r_k + O(r_k^2)\right]\right).\no
\end{align}

We prove first Part (i). Fix $\eta>0.$ Since $\sum_{k=i}^\infty r_k\to 0$ as $i\to\infty,$ there exists $i^*\ge 1$ such that for all $ j\ge t-1\ge i^{*},$
\begin{align*}
    1 - \eta < \exp\left(\sum_{k=t}^j \left[r_k + O\z(r_k^2\y)\right]\right) < 1 + \eta.
\end{align*}
This implies
\begin{align}\label{meta}
    1 - \eta \le m_t \cdots m_j \le 1 + \eta, \quad j \ge t-1 \ge i^*.
\end{align}
Consequently,
\begin{align}\label{rhilu}
    (j - i + 1)(1 - \eta) \le \rho(i,j) \le (j - i + 1)(1 + \eta), \quad  j \ge i \ge i^*.
\end{align}
Choose $N_1 \ge i^*$ such that $\frac{j - i + 1}{j - i} < 1 + \eta$ for all $i \ge N_1$ and $j \ge N_1 + i$. Then,
\begin{align}\label{rije1}
    (1 - \eta)^2 \le \frac{\rho(i,j)}{j - i} \le (1 + \eta)^2, \quad  j \ge i + N_1, \, i \ge N_1.
\end{align}
For $i > i^*$, note that
\begin{align*}
    \rho(0,i) = \sum_{k=1}^{i+1} m_k \cdots m_i = m_{i^*+1} \cdots m_i \sum_{k=1}^{i^*} m_k \cdots m_{i^*} + \rho(i^*,i).
\end{align*}
By \eqref{meta} and \eqref{rhilu},
\begin{align*}
    (C^* + i - i^* + 1)(1 - \eta) \le \rho(0,i) \le (C^* + i - i^* + 1)(1 + \eta), \quad i > i^*,
\end{align*}
where $C^* = \sum_{k=1}^{i^*} m_k \cdots m_{i^*}$. Choose $N_2 > i^*$ such that
\begin{align*}
    1 - \eta \le \frac{C^* + i - i^* + 1}{i} \le 1 + \eta, \quad   i\ge N_2.
\end{align*}
Then,
\begin{align}\label{rie2}
    (1 - \eta)^2 \le \frac{\rho(0,i)}{i} \le (1 + \eta)^2, \quad i \ge N_2.
\end{align}
For any $\varepsilon > 0$, choose $\eta$ sufficiently small so that
\begin{align*}
    1 - \varepsilon < (1 + \eta)^{-2} < (1 - \eta)^{-2} < 1 + \varepsilon,
\end{align*}
and set $N = N_1 \vee N_2$. From \eqref{rije1} and \eqref{rie2}, we obtain \eqref{d0n}. This completes the proof of Part (i).

  Next, we turn to the proof of Part (ii). Fix $B\in [0,1)$ and assume  $r_i=B/i$ for $i\ge1.$ By \eqref{46hyyr4},
  \begin{align}\label{46hyyr42}
    m_t \cdots m_j = \exp\left(\sum_{k=t}^j \left[\frac{B}{k} + O\z(k^{-2}\y)\right]\right), \quad j \ge t-1 \ge 0.
\end{align}
Fix $\eta > 0$. There exists $i^{**} \ge 1$ such that for all $j \ge t-1 \ge i^{**}$,
\begin{align*}
    1 - \eta < \exp\left(\sum_{k=t}^j O\z(k^{-2}\y)\right) < 1 + \eta.
\end{align*}
Observe that for $j \ge t-1 > 0$,
\begin{align*}
    B \log\left(\frac{j+1}{t}\right) \le \sum_{k=t}^j \frac{B}{k} \le B \log\left(\frac{j}{t-1}\right).
\end{align*}
Combining this with \eqref{46hyyr42}, we obtain for $j \ge t-1 \ge i^{**}$:
\begin{align}\label{min}
    (1 - \eta)\left(\frac{j+1}{t}\right)^B \le m_t \cdots m_j \le (1 + \eta)\left(\frac{j}{t-1}\right)^B.
\end{align}
Thus, for $j \ge i \ge i^{**}$,
\begin{align*}
    \rho(i,j) &\le (1 + \eta) \sum_{t=i+1}^{j+1} \left(\frac{j}{t-1}\right)^B
    = (1 + \eta)j^B \sum_{t=i}^j \frac{1}{t^B}
    \le \frac{1 + \eta}{1 - B} j^B \left(j^{1-B} - (i-1)^{1-B}\right).
\end{align*}
For $j > i + 2/\eta$, by the mean value theorem, there exist $\kappa_1 \in (i, i + 2/\eta)$ and $\kappa_2 \in (i-1, i)$ such that
\begin{align*}
    \frac{j^{1-B} - (i-1)^{1-B}}{j^{1-B} - i^{1-B}}
    \le 1 + \frac{i^{1-B} - (i-1)^{1-B}}{(i + 2/\eta)^{1-B} - i^{1-B}}
    = 1 + \left(\frac{\kappa_1}{\kappa_2}\right)^B \frac{\eta}{2}
    \le 1 + \left(\frac{i + 2/\eta}{i - 1}\right)^B \frac{\eta}{2}.
\end{align*}
Choose $\tilde{N}_1 > i^{**} \vee \lfloor \eta/2 + 1 \rfloor$ such that $\left(\frac{i + 2/\eta}{i - 1}\right)^B < 1 + \eta/2$ for all $i\ge\tilde N_1$. Then,
\begin{align*}
    \frac{j^{1-B} - (i-1)^{1-B}}{j^{1-B} - i^{1-B}} \le 1 + \frac{\eta}{2}(1 + \eta/2) < 1 + \eta, \quad  i \ge \tilde{N}_1, \, j - i \ge \tilde{N}_1.
\end{align*}
Consequently,
\begin{align*}
    \rho(i,j) \le \frac{(1 + \eta)^2}{1 - B} j^B \left(j^{1-B} - i^{1-B}\right), \quad  j - i \ge \tilde{N}_1, \, i \ge \tilde{N}_1.
\end{align*}
Similarly, there exists $\tilde{N}_2\ge 1$ such that
\begin{align*}
    \rho(i,j) \ge \frac{(1 - \eta)^2}{1 - B} j^B \left(j^{1-B} - i^{1-B}\right), \quad j - i \ge \tilde{N}_2, \, i \ge \tilde{N}_2.
\end{align*}
Let $\tilde{N}_3 = \tilde{N}_1 \vee \tilde{N}_2$. Then,
\begin{align}\label{dij}
    \frac{(1 - \eta)^2}{1 - B} j^B \left(j^{1-B} - i^{1-B}\right) \le \rho(i,j) \le \frac{(1 + \eta)^2}{1 - B} j^B \left(j^{1-B} - i^{1-B}\right)
\end{align}
for $j - i \ge \tilde{N}_3$ and $i > \tilde{N}_3$. For $\rho(0,i)$, fix $i_0 \ge \tilde{N}_3$ and write
\begin{align*}
    \rho(0,i) = \sum_{k=1}^{i+1} m_k \cdots m_i = m_{i_0+1} \cdots m_i \sum_{k=1}^{i_0} m_k \cdots m_{i_0} + \rho(i_0,i), \quad i> i_0.
\end{align*}
Let $C_0 = \sum_{k=1}^{i_0} m_k \cdots m_{i_0}$. By \eqref{min} and \eqref{dij}, for $i-i_0\ge \tilde N_3$ we have
\begin{align*}
    \rho(0,i) &\le C_0(1 + \eta)\left(\frac{i}{i_0}\right)^B + \frac{(1 + \eta)^2}{1 - B} i^B \left(i^{1-B} - i_0^{1-B}\right), \\
    \rho(0,i) &\ge C_0(1 - \eta)\left(\frac{i+1}{i_0+1}\right)^B + \frac{(1 - \eta)^2}{1 - B} i^B \left(i^{1-B} - i_0^{1-B}\right).
\end{align*}
Therefore, we can choose $N > \tilde{N}_3 + i_0$ such that
\begin{align}\label{pin}
    \frac{i}{1 -B}(1 - \eta)^3 \le \rho(0,i) \le \frac{i}{1 - B} \left(\eta + (1 + \eta)^2\right), \quad i > N.
\end{align}
For any $\varepsilon > 0$, selecting sufficiently small $\eta$ yields \eqref{d1n} from \eqref{dij} and \eqref{pin}.
Part (ii) of the lemma is proved.
  \qed

 We are now ready to complete the proof of Corollary \ref{thz}. To prove Part (i) of Corollary \ref{thz},  for $n\ge 1,$ we set $D(n)=n$ and let  $$S(n):=\sum_{i=1}^n\frac{1}{D(i)}. $$  Clearly, $S(n)\sim \log n$ as $n\rto.$ Thus, the sequence $\{S(n)\}_{n\ge1}$ is  regularly varying  with index $\tau=0.$
By \eqref{pz0}, \eqref{zdij} and Part (i) of Lemma \ref{lmmd}, we may apply Part (ii) of Theorem \ref{thbb} to conclude that
\[\frac{\sum_{j=1}^n\eta_i}{\log n}\overset{d}{\to} \xi\text{ as }n\rto,\] where $\xi\sim \mathrm{Exp}(1).$

Finally, we give the proof of Part (ii) of Corollary \ref{thz}. Utilizing Part (ii) of Lemma \ref{lmmd}, we substitute the parameters $\alpha=1-B$ and $\beta=(1-B)^{-1}$ into Theorem \ref{tha}. This yields
\[\frac{\sum_{j=1}^n\eta_i}{\log n}\overset{d}{\to} \xi\text{ as }n\rto,\]
where $\xi\sim \mathrm{Gamma}(1-B, 1).$
\qed
  \section{Concluding remarks}\label{sc}
    We list here some questions related to Corollaries \ref{thy} and \ref{thz}.

  \begin{itemize}
    \item[(a)]
    Corollary \ref{thy} focuses on critical Galton-Watson processes with geometric offspring distributions. A natural extension, as explored in \cite{w23}, involves considering near-critical branching processes in varying environments. Specifically, let $p_i=\frac{1}{2}+\frac{B}{4 i}$ for $i\ge1$  and define $f_i,\ i\ge1$ as in \eqref{fis}. Let $\{Z_n\}_{n\ge0}$ be a branching process satisfies $Z_0=1$ and $E\z(s^{Z_n}\,\middle|\,Z_0,...,Z_{n-1}\y)=\z(f_i(s)\y)^{Z_{n-1}},\; n\ge 1.$ Try to find the limit distribution of $\#\{1\le t\le n: Z_t=1\}$ (or $\#\{1\le t\le n: Z_t=a\}$ in general for $a\ge 1$) which depends on the value of $B.$

 \item[(b)]
Corollary \ref{thz} addresses the case \( 0 \le B < 1 \). The case \( B \ge 1 \), however, presents distinct behavior. Remark \ref{rmz} indicates that \( C(0) = \{t \ge 1 : Z_t = 0\} \) is almost surely finite, making \( \#C(0) \) a finite random variable. A key challenge is to determine its exact distribution, which likely requires analyzing  the asymptotics of some other multiple sums similar to the one in \eqref{aln}.

    \item[(c)] Both Corollaries \ref{thy} and \ref{thz} assume geometric distributions. As extension, one may generalize these results to arbitrary offspring distributions. A potential strategy involves bounding the probability generating function of \( Z_n \) using linear-fractional functions (see e.g., \cite[Chap. 1]{kv}).
  \end{itemize}

\end{document}